\documentclass[12pt,a4paper,twoside,final,notitlepage, leqno]{article}
\usepackage[english]{babel}
\usepackage[T1]{fontenc}  
\usepackage{graphicx}
\usepackage{amsmath,amsthm,epsfig,amsfonts,bbm}  

\def \smb {{\scriptstyle \bullet }}
\newcommand{\monitem}{ \smallskip \noindent $\bullet$ \quad  } 
\newcommand{\moneq}{\vspace*{-6pt} \begin{equation} \displaystyle } 
\newcommand{\moneqstar}{\vspace*{-6pt} \begin{equation*} \displaystyle } 
\newcommand{\monendstar}{\vspace*{-6pt} \end{equation*}   }
\newcommand{\monend}{\vspace*{-6pt} \end{equation}   }

\def\N{{\rm I}\! {\rm N}}
\def\R{{\rm I}\! {\rm R}}

\def\ib#1{_{_{\scriptstyle{#1}}}}
\def\fT{f\ib{\!\cal{T}}}

\def\hT{h\ib{\cal{T}}}
\def\uT{u\ib{\cal{T}}}
\def\UT{U\ib{\!\cal{T}}}
\def\pT{p\ib{\cal{T}}}
\def\PT{P\ib{\!\cal{T}}}

\def\fT{f\ib{\cal{T}}}
\def\PT{P\ib{\!\!\cal{T}}}
\def\QT{Q\ib{\!\!\cal{T}}}
\def\nT{n\ib{\!\cal{T}}}
\def\xiT{\xi\ib{\!\cal{T}}}
\def\PiT{\Pi\ib{\!\cal{T}}}
\def\MT{M\ib{\!\cal{T}}}
\def\sT{s\ib{\!\cal{T}}}
\def\mT{m\ib{\!\cal{T}}}
\def\wT{w\ib{\!\cal{T}}}
\def\muT{\mu\ib{\!\cal{T}}}

\def\nor#1{\setbox1=\hbox{\kern 3pt{#1}\kern 3pt}%
\dimen1=\ht1 \advance\dimen1 by 0.1pt \dimen2=\dp1 \advance\dimen2 by 0.1pt
\setbox1=\hbox{\kern 1pt  \vrule \kern 2pt \vrule height\dimen1 depth\dimen2\box1
\vrule
\kern 2pt \vrule \kern 1pt  }%
\advance\dimen1 by .1pt \ht1=\dimen1
\advance \dimen2 by .01pt \dp1=\dimen2 \box1 \relax}

\def\mod#1{\setbox1=\hbox{\kern 3pt{#1}\kern 3pt}%
\dimen1=\ht1 \advance\dimen1 by 0.1pt \dimen2=\dp1 \advance\dimen2 by 0.1pt
\setbox1=\hbox{\vrule height\dimen1 depth\dimen2\box1\vrule}%
\advance\dimen1 by .1pt \ht1=\dimen1
\advance \dimen2 by .01pt \dp1=\dimen2 \box1 \relax}

\def\section*#1{}

\def\resume{\if@twocolumn
\section*{R\'esum\'e}
\else \small
\quotation{\bf \it R\'esum\'e \rule[1mm]{1.5mm}{0.2mm}\vspace{0pt}}
\fi}
\def\endresume{\if@twocolumn\else\endquotation\fi}
\def\abstract{\if@twocolumn
\noindent\section*{{\bf Abstract}}
\else \small
\quotation{\noindent \bf {Abstract.} \rule[1mm]{1.5mm}{0.2mm}\vspace{0pt}}
\fi}
\def\endabstract{\if@twocolumn\else\endquotation\fi}

\usepackage{fancyhdr}
\renewcommand{\headrulewidth}{0pt}
\begin{document}

\fancypagestyle{plain}{ \fancyfoot{} \renewcommand{\footrulewidth}{0pt}}
\fancypagestyle{plain}{ \fancyhead{} \renewcommand{\headrulewidth}{0pt}} 

\bibliographystyle{alpha}

\title {\bf \LARGE  Finite volumes and mixed Petrov-Galerkin  \\~ 
 finite elements : the unidimensional problem.    } 

\author { {   \large   Fran\c{c}ois Dubois$^{1,2}$ } \\   
  {\it \small $^1$    Department of Mathematics, University  Paris Sud,}\\ 
  {\it \small   B\^at. 425, F-91405 Orsay Cedex, France. }\\ 
  {\it \small $^2$   Conservatoire National des Arts et M\'etiers, Paris, France,  } \\ 
  {\it \small     Structural Mechanics and Coupled Systems Laboratory.  } \\   
  { \rm  \small francois.dubois@math.u-psud.fr. }   }
 
\date   { {  \rm  04 October 1999} 
\footnote{\rm  \small $\,\,$
Article published in  {\it  Numerical Methods for Partial Differential Equations}, 
volume 16, issue~3, pages 335-360, May 2000. 
Edition 03 January 2014. } }

\maketitle
\renewcommand{\baselinestretch}{1.}

\markboth{Fran\c cois Dubois}
{Finite volumes and mixed Petrov-Galerkin finite elements}

\vspace {-.5 cm} 
\noindent  {\bf Abstract. } \qquad 
For Laplace operator in one space dimension, we propose to formulate the
heuristic finite volume method with the help of mixed Petrov-Galerkin finite
elements. Weighting functions for gradient discretization are parameterized by
some function $\,\psi : [0,1] \rightarrow \R\,$. We propose for this function
$\,\psi\,$ a compatibility interpolation condition and we prove that such a
condition is equivalent to the inf-sup property when studying stability of the
numerical scheme. In the case of stable scheme and under
two distinct hypotheses concerning the regularity of the solution, we
demonstrate convergence of the finite volume method in
appropriate Hilbert spaces and with optimal order of accuracy.

\smallskip  \noindent  {\bf R\'esum\'e. } \qquad 
Dans le cas de l'op\'erateur de Laplace \`a une dimension d'espace,
nous proposons de formuler la m\'ethode heuristique des volumes finis \`a l'aide
d'\'el\'ements finis mixtes dans une variante Petrov-Galerkin o\`u les fonctions
de poids pour la discr\'etisation du gradient sont param\'etr\'ees par une
fonction $\,\psi : [0,1] \rightarrow \R\,$. Nous proposons pour  cette fonction
$\,\psi\,$ une condition de compatibilit\'e d'interpolation qui s'av\`ere
\'equivalente \`a la condition inf-sup pour l'\'etude de la stabilit\'e du
sch\'ema. Dans ce dernier cas et  sous deux hypoth\`eses distinctes concernant la
r\'egularit\'e de la solution, nous d\'emontrons la convergence de la
m\'ethode des volumes finis dans les espaces de Hilbert appropri\'es et
avec un ordre optimal de pr\'ecision.


\smallskip \noindent     {\bf AMS (MOS) classification}: 65N30.

\fancyfoot[C]{\oldstylenums{\thepage}}
\fancyhead[OC]{\sc{Finite volumes and mixed Petrov-Galerkin finite elements}} 









\bigskip \bigskip

\bigskip \bigskip  \newpage \noindent {\bf \large 1) \quad Introduction } 

\monitem 
We study in this paper the approximation of the homogeneous Dirichlet problem for
Poisson equation on the interval $\Omega = \, ]0, 1[$ : 

\smallskip \noindent  (1.1) $\qquad \displaystyle
-\Delta u \quad \equiv \quad {{{\rm d}^{2}u} \over {{\rm d}x^{2}}} \quad = f 
\qquad {\rm in} \,\,    \Omega $


\smallskip \noindent  (1.2) $\qquad \displaystyle
\,\,\, u \qquad = \,\, 0  \qquad  \qquad {\rm on \,\, the \,\, boundary  \,\, \partial
\Omega
\,\, of} \,\, \Omega \,$

\smallskip \noindent 
with the finite volume method. Following, {\it e.g.} Patankar [Pa80], this numerical
method is defined as follows. Consider a ``triangulation'' $\cal{T}$ of the domain
$\Omega$ composed with $(n+1)$ points :

\smallskip \noindent  (1.3) $\qquad \displaystyle
{\cal T} \,\,=\,\, \{ 0 \,=\, x_0 \,{\rm <}\,  x_1 \,{\rm <}\, x_2 \,{\rm <}\, \cdots
\,{\rm <}\,  x_{n-1} \,{\rm <}\, x_n  \,=\, 1\}\,.$

\smallskip \noindent 
The unknowns are the mean values $u_{j+1/2}$ ($j=0, 1, \cdots, n-1$) in each
element $K$ of the mesh $\cal{T}$, with $K$ of the form $K_{j+1/2} \, = \, \, ]x_j,
x_{j+1}[$ :

\smallskip \noindent  (1.4) $\qquad \displaystyle
u_{j+1/2} \quad \approx \quad {{1} \over {x_{j+1} - x_{j}}}  \int_{
\displaystyle x_{j}}^{ \displaystyle x_{j+1}}  u(x) \, {\rm d}x \,. $

\smallskip \noindent 
From these $n$ values, the method proposes an heuristic evaluation of the gradient
$\displaystyle p = {\rm grad} \, u = {{{\rm d}u} \over {{\rm d}x}} \,\, $ at vertex
$x_{j}$ with the help of finite differences :

\smallskip \noindent  (1.5) $\qquad \displaystyle
p_{j} \,\,=\,\,  {{1} \over {h_{j}}} (u_{j+1/2}-u_{j-1/2})\, , \quad j= 0, 1, \cdots, n
\,$

\smallskip \noindent  (1.6) $\qquad \displaystyle
 u_{-1/2} \,\,=\,\, u_{n+1/2}=0  \,$
 
\smallskip \noindent  
to take into account the boundary condition (1.2) ; the length $h_{j+1/2}$ of
interval $]x_j, x_{j+1}[$ is defined by

\smallskip \noindent  (1.7) $\qquad \displaystyle
h_{j+1/2} \,\, = \,\,  x_{j+1} - x_{j}  $


\smallskip \noindent
and distance $h_{j}$ between the centers of two cells $K_{j-1/2}$ and $K_{j+1/2}$
 satisfy the relations

\setbox21=\hbox{$\displaystyle 
h_{0} \,\,=\,\, {{1}\over{2}}  \, h_{1/2}\, $}
\setbox22=\hbox{$\displaystyle \quad 
 h_{j} \,\,=\,\,  {{1}\over{2}}  (h_{j-1/2} + h_{j+1/2})\,,\qquad
j=1,\cdots, n-1 $}
\setbox23=\hbox{$\displaystyle 
 h_{n}  \,\,=\,\, {{1}\over{2}} \, h_{n-1/2}  \,\, . $}
\setbox30= \vbox {\halign{#\cr \box21 \cr   \box22 \cr   \box23 \cr  }}
\setbox31= \hbox{ $\vcenter {\box30} $}
\setbox44=\hbox{\noindent  (1.8) $\displaystyle  \quad   \left\{ \box31 \right. $}  
\smallskip \noindent $ \box44 $

\smallskip \noindent
When $\, p_{j} \, $ is known at vertex $x_{j}$, an integration of the ``conservation
law''\break
 $\displaystyle  \,\, {{\rm div} (p) + f }\equiv {{{ {\rm d}p}\over {{\rm d}x}} + f
}\,=\,  0 \quad $ over the interval $K_{j+1/2}$ takes the following form 

\smallskip \noindent  (1.9) $\qquad \displaystyle
{{1} \over {h_{j+1/2}}} (p_{j+1}-p_{j}) + {{1} \over {h_{j+1/2}}} \int_{ \displaystyle
x_{j}}^{ \displaystyle x_{j+1}}  f(x) \,{\rm d}x \,\,=\,\,  0 \,, \quad j= 0,\cdots,
n-1 
\,$

\smallskip \noindent  
and defines $n$ equations that ``closes'' the problem. This method is very popular,
gives the classical three point finite difference scheme 

\smallskip \noindent  (1.10) $\qquad \displaystyle
{{1} \over {h}} (-u_{j-1/2}+ 2 \, u_{j+1/2}-u_{j+3/2}) \,\,=\,\, {{1} \over {h}}
\int_{\displaystyle x_{j}}^{\displaystyle x_{j+1}}  f(x) \,{\rm d}x  \,, \quad j=
0,\cdots, n-1 \,$

\smallskip \noindent 
for uniform meshes ($h_{j+1/2}\equiv h$ for each $j$), but the numerical 
analysis is difficult in the general case. First tentative was due to Gallou\"et
[Ga92] and weak star topology in space $L^{\infty} (\Omega)$ has been necessary
to take into account the possibility for meshes to ``jump'' abruptly from one value
$ \, h_{j-1/2} \, $ to an other $ \, h_{j+1/2} \, $.

\bigskip \noindent $\bullet \quad$ 
On the other hand, the mixed finite element method proposed by Raviart and Thomas
[RT77] introduces  approximate discrete finite element spaces. Let ${\cal T}$ be a mesh
given at relation (1.3) and $P_{1}$ be the space of polynomials of total degree
$\leq 1$. We set

\smallskip \noindent  (1.11) $\qquad \displaystyle
\UT = \{ u : \Omega \mapsto {\R}, \,\,\forall \,K \in {\cal{T}},\,\,
u_{\mid_{K}} \in {\R}\} \hfill \,$

\smallskip \noindent  (1.12) $\qquad \displaystyle
\PT = \{ p : \overline{\Omega} \mapsto {\R}, \,\, p \,\, {\rm continuous \,\ on} \,\
\Omega,\,\, \forall \, K \in {\cal{T}},\,\, p_{\mid_{K}} \in P_{1}\}\,.$

\smallskip \noindent 
The mixed finite element method consists in solving the problem (1.13)-(1.15) with

\smallskip \noindent  (1.13) $\qquad \displaystyle
\uT \in \UT \, ,\,\, \pT \in \PT \,$

\smallskip \noindent  (1.14) $\qquad \displaystyle
(\pT \,,\, q) + (\uT \,,\,{\rm div}\,q) \,\,=\,\,  0 \,, \qquad \forall \,q \in \PT \,$

\smallskip \noindent  (1.15) $\qquad \displaystyle
({\rm div} \, \pT \,,\, v) + (f \,,\,  v) \,\,\,\,\,=\,\, 0 \,, \qquad \forall \,v \in
\UT \,.$

\smallskip \noindent 
When we explicit the basis  $\chi_{j+1/2} \, (j = 0, 1, \cdots, n-1)$ of linear
space $\UT$ ($\chi_{j+1/2}$ is the numerical function equal to $1$ in $K_{j+1/2}$
and equal to $0$ elsewhere) and the basis $\varphi_{j} \, (j = 0, 1, 2, \cdots,
n)$ of space $\PT$ (recall that $\varphi_{j}$ belongs to space $\PT$ and satisfies
the Kroneker condition  $\varphi_{j}(x_{k}) = \delta_{j,k}$ (for $j$ and $k = 0,
1, 2, \cdots, n$), we introduce vectorial unknowns $\uT$ and $\pT$ according to
the relations

\smallskip \noindent  (1.16) $\qquad \displaystyle
\uT \,\,=\,\,  \sum_{j=0}^{n-1} u_{j+1/2} \, \, \chi_{j+1/2} \,$

\smallskip \noindent  (1.17) $\qquad \displaystyle
\pT =  \sum_{j=o}^{n} p_{j} \,\,\varphi_{j} \,$

\smallskip \noindent
and writing again $\uT$ (respectively $\pT$) the vector in ${\R}^{n}$
(respectively in ${\R}^{n+1}$) composed by the numbers $u_{j+1/2}$ (respectively 
$p_{j}$), system (1.14)-(1.15) takes the form

\smallskip \noindent  (1.18) $\qquad \displaystyle
\left\{ \begin{array} {rcl}  \displaystyle  
M \, \pT & +  \, B^{\displaystyle \rm t} \, \, \uT & = 0  \\ \displaystyle  
B \, \pT   && =  -\fT  \end{array} \right. $

\smallskip \noindent 
with

\smallskip \noindent  (1.19) $\qquad \displaystyle
\fT \,\,=\,\,   \sum_{j=0}^{n-1}  f_{j+1/2} \, \, \chi_{j+1/2} \quad \equiv
\quad \sum_{j=0}^{n-1}  (f,\chi_{j+1/2})  \, \, \chi_{j+1/2} \,. \, $

\smallskip \noindent  
The notations $ \, (\smb ,\smb) \, $ and $\,B^{\displaystyle \rm t} \,$ define
respectively  the scalar product in $L^{2}(\Omega)$ and  the transpose of matrix $B$.
First equation in (1.18) introduces the so-called mass matrix $M$ and  gradient
matrix $B^{\displaystyle \rm t}$ according to formulae

\smallskip \noindent  (1.20) $\qquad 
\left\{ \begin{array} {rll}  \displaystyle  
 M_{j,k}  &\,=\,  (\varphi_{j} \,,\, \varphi_{k}) \,, &  0\leq j \leq n, \quad  0\leq k\leq n
 \\ \displaystyle  
B^{\displaystyle \rm t}_{j,l}  &\,=\, (\chi_{l+1/2}\,,\, {\rm div} \, \varphi_{j}) \,,  
&   0\leq j\leq n, \quad  \leq l\leq n-1 \, . 
 \end{array} \right.  $ 


\smallskip \noindent 
and second equation of (1.18) introduces the {\rm div}ergence matrix $B$ which is the
transpose of the gradient matrix  $ \, B^{\displaystyle \rm t} \, $. The advantage of
mixed formulation is that the numerical analysis is well known [RT77] : the error
$\, \parallel u-\uT\parallel\ib{0} + \parallel p-\pT\parallel\ib{1}\, $ is of order $1$
when the mech size $\hT \equiv \sup_{j}\, h_{j+1/2}$ tends to zero when solution
$u$ of problem (1.1)-(1.2) is sufficiently regular. The main drawback of mixed
finite elements is that system (1.18) is more difficult to solve than system
(1.5)-(1.9) and for this reason, the finite volume method remains very popular.

\bigskip \noindent $\bullet \quad$ 
We focus on the details of non nulls terms of tridiagonal mass matrix ; we have

\smallskip \noindent  (1.21) $\qquad \displaystyle
M_{j,j} \,\, = \,\,  {{2}\over{3}} \, h_{j} \,$

\smallskip \noindent  (1.22) $\qquad \displaystyle
M_{j,j+1} \,\, = \,\, M_{j,j-1} = {{1}\over{6}} \, h_{j+1/2}  \,$

\smallskip \noindent
and therefore

\smallskip \noindent  (1.23) $\qquad \displaystyle
h_{j}  \,\, = \,\,  \sum_{k=0}^{n} M_{j,k}\,,\qquad j=0, 1,\cdots, n \,.$

\smallskip \noindent 
We remark that equation (1.5) is just obtained by the ``mass lumping'' of the first
equation of system (1.18), replacing this equation by the diagonal matrix 
$h_{j}\,\delta_{j,k}$. We refer to Baranger, Ma\^{\i}tre and Oudin [BMO96] for recent
developments of this idea in one and two space dimensions. 

\bigskip \noindent $\bullet \quad$ 
In the following of this article, we show that mixed finite element formulation
(1.13)-(1.15) can be adapted in a Petrov-Galerkin way in order to recover both
simple numerical analysis in classical Hilbert spaces.
Let $d$ be some integer $\geq\,1$ and $\Omega$ be a bounded open set in
${\R}^{d}\,.$ We will denote by $L^{2}(\Omega)$ (or $L^{2}(0,1)$ in one space
dimension when $\Omega$ = ]0,1[) the Hilbert space composed by squarely integrable
functions and by $\parallel\smb\parallel\ib{0}$ the associated norm :

\smallskip \noindent  (1.24) $\qquad \displaystyle
\parallel v \parallel\ib{0}  \, \equiv \, \Bigl( \int_{\Omega} \mid v
\mid^{2} \, {\rm d}x \Bigr) ^{1/2} \quad < \infty  \quad ; $

\smallskip \noindent 
 the scalar product is simply noted with parentheses :

\smallskip \noindent  (1.25) $\qquad \displaystyle
(v,w) \,\,=\,\,  \int_{\Omega}  v(x)\, w(x)\,{\rm d}x \quad .$

\smallskip \noindent 
The Sobolev space  $H^{1}(\Omega)$ is composed with functions in $L^{2}(\Omega)$
whose weak derivatives belong also to space $L^{2}(\Omega)$. The associated norm
is denoted by $\parallel \smb \parallel \ib{1}$ and is defined according to

\smallskip \noindent  (1.26) $\qquad \displaystyle
\parallel v \parallel \ib{1} \, \equiv \, \bigl( \parallel v \parallel\ib{0}^{2}
\, + \, \parallel {\rm grad} \,v \parallel\ib{0}^{2} \bigr) ^{1/2} \,\,, $

\smallskip \noindent
with $\displaystyle \, {\rm grad} \, v = \Bigl( {{\partial v}\over{\partial
x_{1}}},\cdots,{{\partial v}\over{\partial x_{d}}} \Bigr)^{\displaystyle \rm t} \,\,$
and  $\displaystyle \,\parallel {\rm grad} \,v \parallel\ib{0}^{2}\,=\,\sum_{j=1}^{d} 
\parallel {{\partial v}\over{\partial x_{j}}} \parallel\ib{0} ^{2}$. Subspace
$H^{1}_{0}(\Omega)$ of space $H^{1}(\Omega)$ is composed by functions of
$H^{1}(\Omega)$ whose trace values on the boundary $\partial\Omega$ is
identically equal to zero. We will denote by $\,\mid\smb\mid\ib{1}\,$ the
so-called semi-norm associated with space $\,H^{1}_{0}(\Omega)\,$ :  $  \displaystyle
\, \,\mid v \mid\ib{1}^{2}\quad \equiv \quad \parallel {\rm grad} \,v
\parallel\ib{0}^{2}\,. $ The topological dual space of  $\,H^{1}_{0}(\Omega)\,$ is
denoted by $\,H^{-1}(\Omega)\,$ ; note that this space contains $\,L^{2}(\Omega)\,$ but
contains also distributions that can not be represented by functions.

\bigskip \noindent $\bullet \quad$ 
We will use also Sobolev space $\,H^{2}(\Omega)\,$, composed with functions
$\,v\!\in H^{1}(\Omega)\,$ whose gradient also belongs to  $\,H^{1}(\Omega)\,$
 and the associated norm and semi-norm are defined  by the relations

\setbox30=\hbox {$\displaystyle   {{\partial^{2} v}\over{\partial x_{i}
\partial x_{j}}} $}
\smallskip \noindent   $  \displaystyle
\parallel v \parallel \ib{2}^{2}\,\,\equiv \,\, \parallel v \parallel\ib{0}^{2}
\, + \, \parallel {\rm grad} \,v \parallel\ib{1}^{2}\,\,=\,\,
\parallel v \parallel\ib{0}^{2} \, + \, \parallel {\rm grad}
\,v\parallel\ib{0}^{2}\,+\,\sum_{1\leq i,j\leq n}
\nor{\box30} \ib{0}^{2} \,,$
\setbox30=\hbox {$\displaystyle   {{\partial^{2} v}\over{\partial x_{i} \partial
x_{j}}} $}

\smallskip \noindent $ \displaystyle
\mid v \mid \ib{2}\,\,\equiv \,\,\Biggl( \,\sum_{1\leq i,j\leq n}
\nor{\box30} \ib{0}^{2} \, \Biggr) ^{1/2}\,. $

\smallskip \noindent
For mathematical foundation about Sobolev spaces, we refer {\it i.e.} to Adams [Ad75].

\bigskip \noindent $\bullet \quad$ 
The Sobolev space $H({\rm div},\Omega)$ is composed by vector fields
$q = (q_{1},\cdots,q_{d})^{t}\in \bigl(L^{2}(\Omega)\bigr)^{d}$ whose divergence
$\displaystyle \,\,{\rm div}\,q\,\equiv\,\sum_{j=1}^{d}{{\partial q_{j}}\over{\partial
x_{j}}}\,\,$ is  in space $L^{2}(\Omega)$. The norm in space $H({\rm div},\Omega)$ is
denoted by  $\parallel \smb \parallel \ib{{\rm div}}$ and satisfies the natural
relation :

\smallskip \noindent  (1.27) $\qquad \displaystyle
\parallel q \parallel \ib{{\rm div}} \, \equiv \, \Bigl( \,\sum_{j=1}^{d} 
\parallel q_{j} \parallel\ib{0} ^{2} + \parallel {\rm div}\,q \parallel\ib{0}
^{2} \, \Bigr) ^{1/2}\,.$

\smallskip \noindent 
We will often use the product space  $ \, V \, \equiv \, L^{2}(\Omega) \times
H({\rm div},\Omega) \, $ composed by pairs $\eta$ of the form

\smallskip \noindent  (1.28) $\qquad \displaystyle
\eta = (v,q) \quad \in \, L^{2}(\Omega) \times H({\rm div},\Omega)  \,$

\smallskip \noindent
and its natural associated norm satisfies 

\smallskip \noindent  (1.29) $\qquad \displaystyle
\parallel \eta \parallel\ib{V} \, \equiv \, \bigl( \parallel v \parallel\ib{0}^{2}
\, + \, \parallel q \parallel\ib{{\rm div}}^{2} \bigr) ^{1/2} = \bigl( \parallel v
\parallel\ib{0}^{2} \, + \, \parallel q \parallel\ib{0}^{2} + 
\, \parallel {\rm div} \,q \parallel\ib{0}^{2} \bigr) ^{1/2}\,.$

\smallskip \noindent 
without more explicitation.  In one space dimension, the spaces $H({\rm div},]0,1[)$
and  $H^{1}(0,1)$ are identical and we have in this case

\smallskip \noindent  (1.30) $\qquad \displaystyle
\parallel q \parallel\ib{1} \equiv \bigl( \parallel q \parallel\ib{0}^{2} \,
+ \, \parallel {\rm div} \,q \parallel\ib{0}^{2} \bigr) ^{1/2} \,.$

\bigskip \bigskip   \noindent {\bf \large 2) \quad Continuous Petrov-Galerkin formulation} 

\monitem 
We recall in this section the Petrov-Galerkin formulation of problem (1.1)-(1.2)
in the continuous case. Let $d$ be some integer $\geq 1$ and $\Omega \subset
{\R}^{d}$  be a bounded domain with boundary $\partial \Omega$, $u$ be the
solution for the Dirichlet problem for Poisson equation~(2.1)

\setbox21=\hbox{$\displaystyle 
-\sum_{j=1}^{d} \, {{\partial^{2} u}\over{\partial x_{j}^{2}}} \,\,\equiv \,\, 
 -\Delta u \,\,= \,\, f  \qquad {\rm in} \,\,   \Omega, $}
\setbox22=\hbox{$\displaystyle 
\qquad \, u \qquad = \,\, 0 \qquad \qquad \qquad  \,\, {\rm on}  \,\,
\partial \Omega \,\,.  $}
\setbox30= \vbox {\halign{#\cr \box21 \cr   \box22 \cr    }}
\setbox31= \hbox{ $\vcenter {\box30} $}
\setbox44=\hbox{\noindent  (2.1) $\displaystyle  \quad   \left\{ \box31 \right. $}  
\smallskip \noindent $ \box44 $

\smallskip \noindent
First equation of (2.1) can be splitted into two equations of degree $1$ :

\setbox21=\hbox{$\displaystyle 
p \,\,= \,\,  {\rm grad} \, u \qquad \quad \, {\rm in}\quad \Omega  $}
\setbox22=\hbox{$\displaystyle 
{\rm div} \, p + f \,\,= \,\,  0  \qquad {\rm in} \quad  \Omega \,. $}
\setbox30= \vbox {\halign{#\cr \box21 \cr   \box22 \cr    }}
\setbox31= \hbox{ $\vcenter {\box30} $}
\setbox44=\hbox{\noindent  (2.2) $\displaystyle  \quad   \left\{ \box31 \right. $}  
\smallskip \noindent $ \box44 $

\smallskip \noindent
We multiply the first equation of (2.2) by a test function $q \in H({\rm
div},\Omega)$  and second equation of (2.2) by a test function
$v \in L^{2}(\Omega)$. We integrate by parts the right hand side of the first
equation and use the boundary condition in (2.1) to drop out the boundary term.
We  sum the two results and obtain

\smallskip \noindent  (2.3) $\qquad \displaystyle
(u,p) \,  \equiv \, \xi \, \,, \quad  \xi \in \, V \,  \equiv \,  L^{2}(\Omega)\times
H({\rm div},\Omega) \,$

\smallskip \noindent  (2.4) $\qquad \displaystyle
\gamma (\xi,\eta)\quad = \quad <\sigma,\eta>, \qquad \forall \,
\eta\,\equiv \,(v,q)\, \in V \,$

\smallskip \noindent 
with

\smallskip \noindent  (2.5) $\qquad \displaystyle
\gamma \bigl( (u,p)\,,\,(v,q) \bigr) \,\,=\,\, (p\,,\,q) + (u\,,\, {\rm div}\,q) + 
({\rm div}\,p\,,\, v) \,$

\smallskip \noindent  (2.6) $\qquad \displaystyle
<\sigma \,,\, (v,q)>   \,\,=\,\,  -(f\,,\,v) \,.$

\smallskip \noindent
We have the following theorem, due to Babu\u ska [Ba71].

\bigskip
\noindent {\bf  Theorem 1. $\quad$   Continuous mixed formulation.}

\noindent 
Let (V,($\smb$,$\smb$)) be a real Hilbert space, $V'$ its topological dual
space, $\gamma : V \times V \rightarrow\R$ be a continuous bilinear form such that
there exists some $\beta > 0$  satisfying the so-called inf-sup condition :

\smallskip \noindent  (2.7) $\qquad \displaystyle
\inf_{\parallel \xi \parallel \ib{\scriptscriptstyle V} = \,1} \, \,\,
\sup_{\parallel \eta \parallel \ib{\scriptscriptstyle V} \leq \, 1} \,\, 
\gamma (\xi,\eta) \,\, \geq\beta $

\smallskip \noindent 
and a non uniform condition at infinity :

\smallskip \noindent  (2.8) $\qquad \displaystyle
\forall \, \eta \in V, \,\, \bigl( \eta \not = 0 \, \Rightarrow \, \sup_{\xi \in V}
\gamma (\xi,\eta) \,=\,  + \infty \bigr) \,.$

\smallskip \noindent 
Then, for each $\, \sigma \in V'$, the  problem of finding $\xi \in V$
satisfying the relations (2.4) has a unique solution which continuously depends on
$\sigma$ :

\smallskip \noindent  (2.9) $\qquad \displaystyle
\parallel \xi \parallel \ib{V} \quad \leq \quad {{1}\over{\beta}}\,\parallel
\sigma \parallel \ib{V'} \,\, . $

\bigskip \noindent 
The proof of this version of Babu\u ska result can  be found  {\it e.g.} 
 in our report [Du97].

\bigskip \noindent $\bullet \quad $
We show now that choices (2.3) and (2.5) for the Poisson equation leads to a
well-posed problem in the sense of Theorem 1, {\it i.e.} that inf-sup condition
(2.7) and ``infinity condition'' (2.8) are both satisfied.


\bigskip 
  \noindent {\bf Proposition 1. $\quad$   Continuous inf-sup and infinity
conditions.}

\noindent 
Let $V$ be equal to $L^{2}\Omega) \times H({\rm div},\Omega)$ and
$\gamma(\smb,\smb)$ be the bilinear form defined at relation (2.5). Then
$\gamma(\smb,\smb)$ satisfies both inf-sup condition (2.7) and infinity
condition (2.8).

\bigskip \noindent  {\bf Proof of proposition 1.}

\monitem 
We first prove inf-sup condition (2.7). Consider $\,\xi=(u,p)\in V$ with a unity
norm :

\smallskip \noindent  (2.10) $\qquad \displaystyle
\parallel \xi \parallel\ib{V}^{2}  \,\equiv  \,\, \parallel u \parallel\ib{0}^{2} 
+ \parallel p \parallel\ib{0}^{2} + \parallel {\rm div}\,p \parallel\ib{0}^{2} 
\,\,=\, \,1\,. $

\smallskip \noindent
Let $\varphi \in H_{0}^{1}(\Omega)$ be the variational solution of the problem

\setbox21=\hbox{$\displaystyle 
\Delta \varphi  \,\,=\,\, u \qquad {\rm in} \,\,   \Omega \,,  $}
\setbox22=\hbox{$\displaystyle 
\varphi\,\,=\,\, 0  \qquad \,\,\,\, {\rm on} \,\, \partial \Omega \,\,. $}
\setbox30= \vbox {\halign{#\cr \box21 \cr   \box22 \cr    }}
\setbox31= \hbox{ $\vcenter {\box30} $}
\setbox44=\hbox{\noindent  (2.11) $\displaystyle  \quad   \left\{ \box31 \right. $}  

\smallskip \noindent $ \box44 $

\smallskip \noindent
This function $\varphi$ continuously depends on function $u$, {\it i.e.} there
exists some constant $C > 0$ independent of $u$ such that

\smallskip \noindent  (2.12) $\qquad \displaystyle
\parallel \varphi \parallel \ib{1} \, \leq \, C \, \parallel u \parallel
\ib{0} \,.$

\smallskip \noindent
Consider some $\, \beta > 0$ satisfying the inequality

\smallskip \noindent  (2.13) $\qquad \displaystyle
\sqrt{1-\beta-(1+C^{2})(\beta + \sqrt{\beta})^{2}} \,\,\geq \,\,\beta \,. $

\smallskip \noindent 
We verify in the following that we can construct $\eta=(v,q)\in V$ with a norm inferior
or equal to $1$ such that inequality (2.7) holds. We distinguish between three cases,
depending on which term among the three in (2.10) is sufficiently large.

\smallskip  \noindent $\bullet \quad $
If we have

\smallskip \noindent  (2.14) $\qquad \displaystyle
\parallel p \parallel\ib{0}^{2} \,\, \geq \,  \beta\,, $

\smallskip \noindent 
we set $\,\eta\equiv (v,q)$ defined by $\,v=-u\,$ and $\,q=p\,$. We have
clearly, according to (2.5), $\gamma(\xi,\eta)=\, \parallel p \parallel\ib{0}^{2}$ and
inequality (2.7) is a direct consequence of (2.14) in this case. 

\smallskip  \noindent $\bullet \quad $
If inequality (2.14) is in defect and if moreover we have

\smallskip \noindent  (2.15) $\qquad \displaystyle
\parallel u \parallel \ib{0} \,\, \geq \,\, \sqrt{1+C^{2}} \, \bigl(\beta +
\sqrt{\beta} \bigr) \,, $

\smallskip \noindent
we set $\,v=0\,$ and 

\smallskip \noindent   $ \displaystyle
q\,\,=\,\,{{1}\over{ \sqrt{1+C^{2}} \,\,\parallel u \parallel \ib{0}}}\,\,
{\rm grad}\,\varphi\, $

\smallskip \noindent
with $\varphi$ introduced in (2.11). Then it follows from relation (2.12) that the
norm $\, \parallel \eta \parallel\ib{V}\,$ of $\eta=(v,q)$ is not greater than~$1$
because $ \displaystyle \parallel q \parallel \ib{0}\,\, \leq\, {{C}\over
{\sqrt{1+C^{2}}}}\,$.  We have moreover 

\smallskip \noindent   $ \displaystyle
\gamma(\xi \,,\, \eta) \,\, \geq \,\, (p \,,\, q) + (u \,,\, {\rm div}\,q) + ({\rm
div}\,p \,,\, v) \,$

\smallskip \noindent   $ \qquad  \quad  \,\,\,\, \displaystyle
\geq \,\, (u\,,\,{\rm div}\,q) \,\,- \,\parallel p \parallel\ib{0} \,\parallel q
\parallel\ib{0} $

\smallskip \noindent   $ \qquad   \quad  \,\,\,\,  \displaystyle
\geq \,\,  {{\parallel u \parallel\ib{0}}\over{\sqrt{1+C^{2}}}} \, \,-\,\,\sqrt{\beta}
\,$

\smallskip \noindent 
and due to (2.15) this last quantity is greater than $\beta$ ; inequality (2.7)
is established in this second case.

\smallskip  \noindent $\bullet \quad $
If inequalities (2.14) and (2.15) are both in defect, we set $\displaystyle \,v={{{\rm
div}\,p}\over{\parallel {\rm div}\,p \parallel \ib{0}}}\,$ and $\,q=0$. Then
$\eta=(v,q)$ is of unity norm and $\gamma(\xi,\eta)=\,\parallel {\rm div}\,p\parallel
\ib{0}$. But from equality (2.10) we have also

\smallskip \noindent   $ \displaystyle
 \parallel {\rm div}\,p\parallel \ib{0}^{2} \,\,=\,\, 1 \,- \parallel u \parallel 
 \ib{0}^{2} - \parallel p \parallel \ib{0}^{2} \, $

\smallskip \noindent   $ \qquad   \qquad  \, \, \displaystyle
\geq \,\, 1 - (1+C^{2})\,(\beta + \sqrt{\beta})^{2}\,-\,\beta \,\, \geq \,\, 
\beta^{2} $ 

\smallskip \noindent 
due to relation (2.13). Then the inf-sup inequality (2.7) is established. 

\smallskip  \noindent $\bullet \quad $
We prove now the infinity condition (2.8). Let $\eta=(v,q)$ be a non-zero pair
of functions  in the product space $L^{2}(\Omega)\times H({\rm div},\Omega)$. We
again distinguish between  three cases.

\smallskip \noindent 
(i) $\quad$ If ${\rm div}\,q\ne0$, we set $u=\lambda\,\, {\rm div}\,q$, $p=0$ and
$\xi=(u,p)$. Then $\,\gamma(\xi,\eta)=$ $=\lambda\,\parallel {\rm div}\,p
\parallel\ib{0}^{2} \, $ tends to $+\infty$ as $\lambda$ tends to  $+\infty$.

\smallskip \noindent 
(ii) $\quad$  If $\,{\rm div}\,q=0$ and $v \ne 0$, let $ \varphi \in
H_{0}^{1}(\Omega)$ be the variational solution of the problem

\smallskip \noindent   $ \displaystyle
\left\{ \begin{array} {rll}  \displaystyle  
\Delta \varphi \, &\,=\,  v   & {\rm in} \,\,   \Omega
 \\ \displaystyle  
\varphi  &\,=\,  0  &{\rm on  \,\, \partial \Omega \,\,}
 \end{array} \right.  $ 
%


\smallskip \noindent 
and $\, \widetilde{p}={\rm grad}\,\varphi$. Then
$\, (\widetilde{p},q)=({\rm grad}\,\varphi,q) \,=\, -(\varphi,{\rm div}\,q)=0$. We set
$u=0$, $p=\lambda\,\widetilde{p}$ and  $\xi=(u,p)$. We have $\gamma(\xi,\eta)=\lambda
({\rm div}\,\widetilde{p} \,,\, v)=\lambda\,\parallel v\parallel \ib{0}^{2}$ which
tends to $+\infty$ as $\lambda$ tends to  $+\infty$.

\smallskip \noindent 
(iii) $\quad$  If ${\rm div}\,q=0$ and $v=0$, vector $q$ is non null by hypothesis. 
Then $u=0$, $p=\lambda \,q$ and  $\xi=(u,p)$ show that $\gamma(\xi,\eta)=(p,q) \,=\,
\lambda\,\parallel q\parallel \ib{0}^{2}$ which tends to $+\infty$ as $\lambda$ tends
to  $+\infty$. Inequality (2.8) is established and the proof of Proposition~1 is
completed. $\hfill  \square \kern0.1mm$

\bigskip \bigskip   \noindent {\bf \large 3) \quad Discrete mixed Petrov-Galerkin 
formulation  for finite volumes} 

\monitem 
We consider again the unidimensional problem (1.1)-(1.2)  on
domain $\Omega=]0,1[$, the mesh ${ \cal T}$ introduced in (1.3), a discrete
approximation space $\,\UT$ of Hilbert space $L^{2}(\Omega)$  defined in (1.11) 
and a discrete finite dimensional approximation space $\PT$ of Sobolev space
$H({\rm div},\Omega)$ defined at relation (1.12). We modify in the following the
mixed finite element formulation (1.13)-(1.15) of problem (1.1)(1.2) and consider
the discrete mixed Petrov-Galerkin formulation : 

\smallskip \noindent  (3.1) $\qquad \displaystyle
\uT \in \UT \, ,\,\, \pT \in \PT  \,$

\smallskip \noindent  (3.2) $\qquad \displaystyle
(\pT \,,\, q) + (\uT\,,\, {\rm div}\,q) = 0 \,, \qquad \forall \,q \in \QT^{\psi} \,$

\smallskip \noindent  (3.3) $\qquad \displaystyle
({\rm div} \, \pT,v) + (f\,,\, v) = 0 \,, \qquad \,\, \forall \,v \in \UT \,.$

\smallskip \noindent
We remark that the only difference with (1.13)-(1.15) consists in the choice of
test function $q$ in relation (3.2) : in the classical mixed formulation, $q$
belongs to space $\PT$ (see relation (1.14)) whereas in the present one, we
suppose in equation (3.2) that $q$ belongs to space $\QT^{\psi}$. The trial
functions (space $\PT$) and the weighting functions (space $\QT^{\psi}$)  for the
discretization of the eqation $\,p={\rm grad}\,u\, $ are now not identical.
Therefore we have replaced a classical mixed formulation by a Petrov-Galerkin
one, in a way suggested several years ago by Hughes [Hu78] and Johnson-N\"avert 
[JN81] for advection-diffusion problems, more recently in a similar context by
Thomas and Trujillo [TT99]. 

\smallskip \noindent
We define the space $\QT^{\psi}$ in the way described below.

\bigskip
\noindent {\bf Definition 1. $\quad$  Space of weighting functions.}

\noindent 
Let $\psi : [0,1] \rightarrow \R\quad$ be a continuous function satisfying the
localization condition

\smallskip \noindent  (3.4) $\qquad \displaystyle
\psi(0) = 0\,,\quad \psi(1)=1 \,, $

\smallskip \noindent 
let $\cal{T}$ be a mesh  given in relation (1.3) and defined by 
vertices $x_{j}$ and finite elements $K$ of the form  $K_{j+1/2} = \,\, ]x_j,
x_{j+1}[$. We define a basis function $\psi_{j}$ of space $\QT^{\psi}$ by affine
transformation of function $\psi$ : 

\setbox21=\hbox{$\displaystyle 
\psi\Bigl( {{x-x_{j-1}}\over{h_{j-1/2}}} \Bigr) \quad  {\rm if}\,\, x_{j-1} \leq
x \leq x_{j} \,$}
\setbox22=\hbox{$\displaystyle 
\psi\Bigl( {{x_{j+1}-x}\over{h_{j+1/2}}} \Bigr) \quad {\rm if}\,\, x_{j} \leq x \leq
x_{j+1} \,$}
\setbox23=\hbox{$\displaystyle 
0 \qquad \qquad \qquad  \,  {\rm elsewhere }\,. $}
\setbox30= \vbox {\halign{#\cr \box21 \cr   \box22 \cr   \box23 \cr  }}
\setbox31= \hbox{ $\vcenter {\box30} $}

\setbox44=\hbox{\noindent  (3.5) $\displaystyle  \quad   \psi_{j}(x) \,\, =
\,\,  \left\{ \box31 \right. $}  

\smallskip \noindent $ \box44 $

\smallskip \noindent 
The space $\QT^{\psi}$ is defined as the set of linear combinaisons of
functions $\psi_{j}$ :

\smallskip \noindent  (3.6) $\qquad \displaystyle
q \in \QT^{\psi}\quad {\rm iff} \quad  \exists \,q_{0},\cdots,q_{n} \in
\R \quad {\rm such \,\,that} \quad q={\sum_{j=0}^{n}} q_{j} \,\psi_{j}
\,\,.  $

\bigskip  \noindent $\bullet \quad$ 
The interest of such weighting functions is to be able to diagonalize the mass
matrix $(\varphi_{i},\psi_{j})$ $(0\leq i,j \leq n)$ composed with the basis
$(\varphi_{i})\ib{0\leq i\leq n}$ of space $\PT$ and the basis
$(\psi_{j})\ib{0\leq j\leq n}$ of linear space  $\QT^{\psi}$. We have the
following result : 

\bigskip 
\noindent {\bf Proposition 2. $\quad$  Orthogonality.}

\noindent 
Let $\psi$ be defined as in definition 1 and satisfying moreover the
orthogonality condition

\smallskip \noindent  (3.7) $\qquad \displaystyle
\int_{0}^{1} (1-x)\,\psi(x)\, {\rm d}x \quad=\quad 0\,. $

\smallskip \noindent
Then the mass matrix $(\varphi_{i},\psi_{j})$ $(0\leq i,j \leq n)$ associated
with equation (3.2) is diagonal :

\smallskip \noindent  (3.8) $\qquad \displaystyle
\exists\,H_{j}\,\in\R\,,\quad (\varphi_{i},\psi_{j}) \,=\, H_{j} \, 
\, \delta_{i,j}\,,\qquad 0\leq i,\, j \leq n\,.$

\bigskip \noindent  {\bf Proof of proposition 2.}

\monitem 
The proof of relation (3.8) is elementary. If $i$ and $j$ are two
different integers, the support of function $\,\varphi_{i}\,\psi_{j}\,$ is
reduced to a null Lebesgue measure set except if $i=j-1$ or $i=j+1$. In the first
case, we have

\smallskip \noindent  $ \displaystyle
\int_{ 0}^{ 1} \varphi_{j-1}(x) \, \psi_{j}(x)\, {\rm
d}x  \,\,= \,\,   \int_{\displaystyle x_{j-1}}^{\displaystyle x_{j}}
\varphi_{j-1}(x) \, \psi_{j}(x)\,  {\rm d}x$

\smallskip \noindent  $ \qquad \qquad  \qquad \qquad  \quad \,\,\, \displaystyle
= \,\,    h_{j-1/2}\,\,  \int_{0}^{1} (1-y)\,\psi(y)\,
{\rm d}y\,\,, \qquad  j=1,\cdots,n \,$

\smallskip \noindent
with the change of variable $\,\,x=x_{j-1}+h_{j-1/2}\,y\,\,$ compatible with
relations (3.5). The last expression in the previous computation is null due to
(3.7).

\smallskip  \noindent $\bullet\quad$ 
In a similar way, in the second case, we have :

\smallskip \noindent  $ \displaystyle
\int_{0}^{1}   \varphi_{j+1}(x) \, \psi_{j}(x)\,{\rm
d}x \,\,= \,\,   \int_{\displaystyle x_{j}}^{\displaystyle x_{j+1}} \,
\varphi_{j+1}(x) \, \psi_{j}(x)\, {\rm d}x \, $

\smallskip \noindent  $ \qquad \qquad  \qquad \qquad  \quad \,\,\, \displaystyle
= \,\,    h_{j+1/2}\,\,  \int_{0}^{1} (1-y)\,\psi(y)\, {\rm d}y\,\,, \qquad 
j=1,\cdots,n \,$

\smallskip \noindent 
with a new variable $y$ defined by the relation $\,\,x=x_{j+1}-h_{j+1/2}\,y\,\,$
and thanks to relation (3.5). The resulting integral remains equal to zero due to
the orthogonality condition (3.7).

\smallskip  \noindent $\bullet\quad$ 
When $j=i$, previous calculations show that

\smallskip \noindent  $ \displaystyle
\int_{0}^{1}   \varphi_{j}(x) \, \psi_{j}(x)\,{\rm d}x \,\,= \,\,   
\int_{\displaystyle x_{j-1}}^{\displaystyle x_{j}}
\varphi_{j}(x) \, \psi_{j}(x)\, {\rm d}x \,\,+\,\, 
\int_{\displaystyle x_{j}}^{\displaystyle x_{j+1}}
\varphi_{j}(x) \, \psi_{j}(x)\, {\rm d}x \,$

\smallskip \noindent  $ \qquad \qquad  \qquad \qquad  \quad \,\,\, \displaystyle
= \,\, \big(h_{j-1/2}\,+\,h_{j+1/2}\big) \, 
\int_{0}^{1} y\,\psi(y)\, {\rm d}y  \,\,,\quad
j=1,\cdots,n-1\,.$

\smallskip \noindent 
If $h_{j}$ is the expression defined in (1.8),  the value of $H_{j}$ is simply
expressed by :

\smallskip \noindent  (3.9) $\qquad \displaystyle
H_{j} \,\,=\,\,  2 \, h_{j}\,\int_{0}^{1} x\,\psi(x)\, {\rm d}x \,\,,\qquad
j=0,\cdots,n  \,$

\smallskip \noindent 
and Proposition 2 is then proven.  $\hfill  \square \kern0.1mm$

\bigskip \bigskip \noindent $\bullet \quad$
We can now specify a choice of shape function $\psi$ in order to recover finite
volumes with mixed Petrov-Galerkin formulation : since relation (3.2) used with
test function $q=\psi_{j}$ shows  (with notations given at relations (1.16) and
(1.17)) :

\smallskip \noindent  (3.10) $\qquad \displaystyle
H_{j}\,p_{j} \,\,=\,\, u_{j+1/2} - u_{j-1/2}\,,\qquad j=0,\cdots,n \,,$

\smallskip \noindent 
the finite volumes are reconstructed if  relation (3.10) is identical to
the heuristic definition (1.5), {\it i.e.} due to (3.9), if we have the following
compatibility condition between finite volumes and mixed Petrov-Galerkin
formulation :

\smallskip \noindent  (3.11) $\qquad \displaystyle
\int_{0}^{1} x\,\psi(x)\, {\rm d}x \,\,=\,\, {1\over2}\,. $

\smallskip \noindent 
The next proposition show that cubic spline function can be choosen as
localization $\psi$ function. 

\bigskip
\noindent {\bf Proposition 3. $\quad$  Spline example.}

\noindent 
Let  $\psi : [0,1] \rightarrow \R\quad$ be a continuous function satisfying the
localization condition (3.4), orthogonality condition (3.7) and the
compatibility condition with finite volumes (3.11). Then function $\psi$ is
uniquely defined if we suppose moreover that  $\psi$ is polynomial of degree
$\leq 3$. We have

\smallskip \noindent  (3.12) $\qquad \displaystyle
\psi(x)\,\,=\,\, {1\over2} + 3\,(2x-1) - {5\over2}\,(2x-1)^{3} \,\,=\,\,
-9x+30x^{2}-20x^{3} \,.$

\bigskip \noindent  {\bf Proof of proposition 3.}

\monitem 
It  is an elementary calculus. First, due to (3.4), it is natural to search
$\psi$ of the form $\,\, \psi(x)=x\,\bigl(1+\alpha(1-x)+\beta(1-x)^{2} \bigr)$.
Secondly it comes simply from (3.7) and (3.11) that 

\smallskip \noindent  $ \displaystyle
\int_{0}^{1} \psi(x)\, {\rm d}x \,\,=\,\,\int_{0}^{1}
x\,\psi(x)\, {\rm d}x \,\,=\,\, {1\over2}\,. $

\smallskip \noindent
Then due to the explicit value of some polynomial integrals

\smallskip \noindent  $ \displaystyle
\int_{0}^{1} x(1-x)\,{\rm d}x = {1\over 6}\,\,,\quad
\int_{0}^{1} x^{2}(1-x)\,{\rm d}x = {1\over
12}\,\,,\quad\int_{0}^{1} x^{2}(1-x)^{2}\,{\rm d}x =
{1\over 30} \,, $

\smallskip \noindent
we can express $\, \displaystyle \int_{0}^{1}
\psi(x)\, {\rm d}x\,\,$ and $\, \displaystyle \int_{0}^{1} x\,\psi(x)\, {\rm
d}x\,\,$ in terms of unknowns $\alpha$ and $\beta$ :

\smallskip \noindent  $ \displaystyle
\int_{0}^{1} \psi(x)\, {\rm d}x= {1\over 2} + {{\alpha}\over 6} + {{\beta}\over
12} \,,\quad \int_{0}^{1} x\,\psi(x)\, {\rm d}x= {1\over 3} + {{\alpha}\over 12}
+ {{\beta}\over 30} \,.  $

\smallskip \noindent
We deduce that $\alpha=10\,$, $\beta=-20\,$ and relation (3.12) holds.  $\hfill 
\square \kern0.1mm$

\bigskip \bigskip   \noindent {\bf \large 4) \quad  Discrete inf-sup condition } 

\monitem  
For unidimensional Poisson equation with homogeneous boundary condition, the
finite volume method is now formulated as a discrete approximation (3.1)-(3.3)
associated with the biliear form $\gamma(\smb,\smb)$ defined in relation
(2.5) and the following finite dimensional subspaces $V_{1}$ and $V_{2}$ of
continuous space $V = L^{2}(\Omega) \times H({\rm div},\Omega)$~:

\smallskip \noindent  (4.1) $\qquad \displaystyle
V_{1} \,\,=  \,\, \UT \times  \PT \,$

\smallskip \noindent  (4.2) $\qquad \displaystyle
V_{2}  \,\,= \,\,  \UT \times   \QT^{\psi} \,\,. $

\smallskip \noindent 
With these notations, problem (3.1) (3.2) (3.3) can be formulated as follows :

\smallskip \noindent  (4.3) $\qquad \displaystyle
\xi_{1} = (\uT,\pT) \,\, \in  \,\,V_{1} \,$

\smallskip \noindent  (4.4) $\qquad \displaystyle
\gamma (\xi_{1} \,,\, \eta)  \,\,= \,\, < \sigma\,,\,  \eta >\,\,,\quad \forall\,
\eta \in V_{2} \,$

\smallskip \noindent 
with linear form $\sigma$ defined in (2.6). We have the following approximation
theorem [Ba71].

\bigskip 
\noindent {\bf Theorem 2. $\quad$  General approximation result.}

\noindent 
Let  $V$ be a real Hilbert space and  $\gamma$ be a continuous bilinear form
like in Theorem~1 with a continuity modulus denoted by M : 

\smallskip \noindent  (4.5) $\qquad \displaystyle
\mid \gamma(\xi,\eta) \mid \,\,\leq\,\, M \, \parallel \xi \parallel\ib{V}\, 
\parallel \eta \parallel\ib{V}\,,\qquad \forall \, \xi,\eta\in V\,. $ 

\smallskip \noindent 
Let $V_{1}$ and $V_{2}$ be two closed subspaces of space $V$ such that we have
the following two properties : on one hand, there exits some constant $\delta$
associated with the  uniform  discrete inf-sup condition

\smallskip \noindent  (4.6) $\qquad \displaystyle
\inf_{\xi\in V_{1},\,\parallel \xi \parallel \ib{\scriptscriptstyle V} = \,1} \,
\,\, \sup_{\eta\in V_{2},\,\parallel \eta \parallel \ib{\scriptscriptstyle V}
\leq \, 1} \,\, \gamma (\xi,\eta) \,\, \geq \, \delta \,$

\smallskip \noindent
and on the other hand, the discrete infinity condition

\smallskip \noindent  (4.7) $\qquad \displaystyle
\forall \,\eta \in V_{2} \setminus \{0\}\,,\quad \sup_{\xi \in
V_{1}} \,\gamma(\xi,\eta) \, = +\infty \,$

\smallskip \noindent
is satisfied. Then problem (4.3)(4.4) has a unique solution
$\,\xi_{1}\!\!\in\!\!V_{1}$. If $\xi$ is the solution of continuous problem
(2.3)(2.4) (obtained simply with  $V_{1}=V_{2}=V$), we have the following control
of the approximation error by the interpolation error :

\smallskip \noindent  (4.8) $\qquad \displaystyle
\parallel \xi - \xi_{1} \parallel\ib{V} \,\,\,\leq\,\,\, \Big( 1+{M\over{\delta}}
\Bigr) \parallel \xi - \zeta \parallel\ib{V} \,,\qquad \forall \, \zeta \in
V_{1}\,. $

\monitem 
Theorem 2 plays an analogous role than the so-called  Cea lemma [Ce64] in
classical analysis of the error for conforming finite elements (Ciarlet-Raviart
[CR72]). It states that when constant
$\delta$ in estimate (4.6) is independent of the choice of spaces $V_{1}$ and
$V_{2}$ (uniform inf-sup discrete condition) the error $\parallel
\xi-\xi_{1}\parallel\ib{V}$ is dominated by the interpolation error
$\inf_{\zeta\in V_{1}} \parallel \xi-\zeta\parallel\ib{V}$, that establishes
convergence with an optimal order when $V_{1}$ is growing more and more towards
space $V$. The two next propositions compare discrete $L^{2}$ norms when
interpolation function $\psi$, satisfying the two conditions (3.4) and (3.7), is
moreover submitted to the following {\bf compatibility interpolation condition}

\smallskip \noindent  (4.9) $\qquad \displaystyle
\psi(\theta) + \psi(1-\theta) \,\equiv\, 1 \,,\qquad \forall \,\theta
\in [0,1] \,$

\smallskip \noindent 
does not satisfy it. Note that for the spline example (3.12), compatibility
interpolation condition was satisfied. We suppose also that the mesh $\cal T$ can
be chosen in the class $\,\cal U_{\alpha,\beta}$ of uniformly regular meshes.

\bigskip
\noindent {\bf Definition 2. $\quad$  Uniformly regular meshes.}

\noindent 
Let $\alpha,\beta$ be two real numbers such that

\smallskip \noindent  (4.10) $\qquad \displaystyle
0 \, < \, \alpha  \, < \, 1  \, < \, \beta\,.$

\smallskip \noindent
The class $\cal U_{\alpha,\beta}$ of uniformly regular meshes is composed by all
the meshes $\cal T$ associated with $\nT$ ($\nT \! \in \N$)  vertices 
$x_{j}^{\cal T}$ satisfying

\smallskip \noindent  (4.11) $\qquad \displaystyle
0=x_{0}^{\cal T}<x_{1}^{\cal T}	<\cdots < x_{\nT-1}^{\cal T} <  
x_{\nT}^{\cal T} = 1 \,$

\smallskip \noindent 
and such that the corresponding measures $h_{j+1/2}^{\cal T}$ of elements
$K_{j+1/2}^{\cal T}$

\smallskip \noindent   $ \displaystyle
h_{j+1/2}^{\cal T}\,\,=\,\, x_{j+1}^{\cal T} - x_{j}^{\cal T}\,\,,\quad
j=0,1\cdots, \nT-1 \,$

\smallskip \noindent 
satisfy the condition

\smallskip \noindent  (4.12) $\qquad \displaystyle
{{\alpha}\over{\nT}} \,\,\leq\,\, h_{j+1/2}^{\cal T} \,\,\leq\,\, {{\beta}\over
{\nT}} \,,\quad \forall \,j=0,1\cdots, \nT-1\,,\quad \forall \, {\cal T} \,\in
\cal U_{\alpha,\beta} \,.$

\bigskip \noindent $\bullet \quad$
We remark that the ratio $\,{h_{j+1/2}^{\cal T}}\,$/$\,{h_{j-1/2}^{\cal T}}\,$
of successive cells has not to be close to $1$ but remains bounded from below by
$ \, \alpha \, $/$ \, \beta\, $ and from above by $ \, \beta \, $/$ \, \alpha
$. We will denote by $\hT$ the maximal stepsize of mesh $\cal T$ : 

\smallskip \noindent  (4.13) $\qquad \displaystyle
\hT = \max_{j=0,1,\cdots,\nT\!-1} h_{j+1/2}^{\cal T}\,. $

\bigskip  \noindent {\bf Proposition 4. $\quad$  Stability when changing the interpolant
function.}

\noindent 
Let $\psi$ be a continuous function $\,[0,1] \rightarrow \R\quad$ satisfying the
conditions (3.4), (3.7) and the compatibility interpolation condition (4.9).
Let $\cal T$ be some mesh of the interval $[0,1]$ composed with $\nT\!=n$
elements, $\PT$ be the space of continuous $P_{1}$ functions associated with mesh
$\cal T$ and defined in (1.12) and $\QT^{\psi}$ be the analogous space, but
associated with the use of $\psi$ for interpolation and defined in (3.5)(3.6). 
Consider $(q_{0}, q_{1},\cdots,q_{n})\,\in \R^{n+1}$, 

\smallskip \noindent  $ \displaystyle
q={\sum_{j=0}^{n}} q_{j} \,\psi_{j} \,\,\in \QT^{\psi}\quad{\rm and}\quad
\widetilde q ={\sum_{j=0}^{n}} q_{j} \,\varphi_{j} \,\, \in \PT\,. $

\smallskip \noindent 
We have the estimations

\smallskip \noindent  (4.14) $\qquad \displaystyle
{3\over2}\,\delta\, \parallel \widetilde q \parallel \ib{0}^{2}\quad \leq
\quad \parallel  q \parallel \ib{0}^{2} \quad \leq \quad
12 \,\widetilde{\delta}\, \parallel \widetilde q \parallel \ib{0}^{2} \,$

\smallskip \noindent 
with strictly positive constants $\delta$ and $\widetilde{\delta}$ defined by

\smallskip \noindent  (4.15) $\qquad \displaystyle
\setbox20= \hbox{$\displaystyle \int_{0}^{1} \psi(\theta)\,\psi(1-\theta)\,{\rm
d}\theta$} \delta \,\,=\,\,  \int_{0}^{1} \bigl(\psi(\theta)\bigr)^{2}\,{\rm d}
\theta  \,\,-\,\, \mod{\box20} \,$

\smallskip \noindent  (4.16) $\qquad \displaystyle
\widetilde{\delta} \quad=\quad \int_{0}^{1} \bigl(\psi(\theta)\bigr)^{2} \, {\rm
d} \theta \,. $

\bigskip \noindent  {\bf Proof of proposition 4.}

\smallskip \noindent $\bullet\quad$ 
It is not immediate that $\delta$ is strictly positive. From Cauchy-Schwarz
inequality we have

\setbox20= \hbox{$\displaystyle \int_{0}^{1}
\psi(\theta)\,\psi(1-\theta)\,{\rm d} \theta$}
\smallskip \noindent  (4.17) $\qquad \displaystyle
\mod{\box20} \quad \leq \quad \int_{0}^{1} \bigl(\psi(\theta)\bigr)^{2} \,
{\rm d} \theta \,$

\smallskip \noindent
which proves that $\delta \geq 0$. If there is exact equality in inequality
(4.17), the case of equality in Cauchy-Schwarz inequality show that the two
functions in the scalar product at the left hand side of (4.17) are proportional :

\smallskip \noindent $ \displaystyle
\exists \,\lambda,\mu \in \R, \,(\lambda,\mu) \ne (0,0)\,,\quad \forall \,
\theta \in [0,1]\,, \quad \lambda \,\psi(1-\theta) + \mu \,\psi(\theta) =
0\,. \,$

\smallskip \noindent
Taking $\theta=0$ in previous inequality, localization condition (3.4) shows that
$\lambda=0$. In a similar manner,  the choice of the particular value
$\,\theta=1\, $  implies $\,\mu=0\,$, which is finally not possible because
$\,(\lambda,\mu) \ne (0,0)\,$. Therefore the equality case in (4.17) is excluded
and $\,\delta > 0\,$. 

\smallskip \noindent $\bullet \quad$
We evaluate now the $L^{2}$ norm of $q={\sum_{j=0}^{n}} \, q_{j} \,\psi_{j}$. We get

\setbox20=\hbox{$\parallel q \parallel \ib{0}^{2}$}
\setbox21=\hbox{$\displaystyle  \sum_{j=0}^{n-1} \, \int_{\displaystyle x_{j}}
^{\displaystyle x_{j+1}} \Bigl(  \, q_{j} \,\, \psi\Bigl({{x_{j+1}-x}
\over{h_{j+1/2}}}\Bigr) +   q_{j+1} \,\, \psi
\Bigl({{x-x_{j}}\over{h_{j+1/2}}}\Bigr) \, \Bigr)^{2} \,{\rm d}x$}
\setbox22=\hbox{$\displaystyle  \sum_{j=0}^{n-1} h_{j+1/2}  \int_{0}^{1} \bigl[ 
q_{j} \, \psi(1\!-\!\theta) \,+\,  q_{j+1} \, \psi(\theta) \bigr] ^{2}\,
{\rm d} \theta\,\,,\quad x=x_{j}+\theta  \, h_{j+1/2} $}
\setbox23=\hbox{$\displaystyle  \sum_{j=0}^{n-1} h_{j+1/2}  
\biggl[ (q_{j}^{2}+q_{j+1}^{2})\,\int_{0}^{1} (\psi(\theta))^{2}\,{\rm d}
\theta \,+\, 2\,q_{j}\,q_{j+1} \int_{0}^{1} \psi(\theta) \,\psi(1\!-\!\theta)
\,{\rm d} \theta \bigg]$}
\setbox25=\hbox{$\displaystyle \int_{0}^{1} \psi(\theta) \,\psi(1\!-\!\theta)
\,{\rm d}\theta$}
\setbox24=\hbox{$\displaystyle  \sum_{j=0}^{n-1} h_{j+1/2} \,\,
(q_{j}^{2}+q_{j+1}^{2}) \,\, \biggl[ \int_{0}^{1} (\psi(\theta))^{2}\,{\rm
d}\theta \, -\, \mod{\box25}\, \biggr] $}
\setbox25=\hbox{$\displaystyle \int_{0}^{1} \psi(\theta) \,\psi(1\!-\!\theta)
\,{\rm d}\theta$}
\setbox26=\hbox{$\displaystyle \qquad \quad + \quad \sum_{j=0}^{n-1} h_{j+1/2}
\,\, \bigl(\mid q_{j} \mid - \mid q_{j+1} \mid \bigr)^{2}\,\, \mod{\box25}$}
\setbox27=\vbox{\halign{#&#&#\cr \box20&$\quad=\quad$&\box21 \cr 
&$\quad=\quad$& \box22 \cr &$\quad=\quad$& \box23 \cr 
&$\quad \geq \quad$& \box24 \cr &$\,$& \box26 \cr}}
\smallskip \noindent $\box27$

\smallskip \noindent  (4.18) $\qquad \displaystyle
\parallel q \parallel \ib{0}^{2} \quad \geq \quad
\delta \,\,  \sum_{j=0}^{n-1} h_{j+1/2} \,\, (q_{j}^{2}+q_{j+1}^{2})\,.$

\smallskip \noindent
We have an analogous inequality concerning $\,\,\widetilde q ={\sum_{j=0}^{n}}
q_{j}\,\varphi_{j}\,$, by replacing the number $\delta$ by its precise value when
$\psi(\smb)$ is replaced by an affine interpolation between data, ie function
$\, \R \ni \theta \mapsto \theta \in \R \,$. We deduce from (4.18) in this
particular case :

\smallskip \noindent  (4.19) $\qquad \displaystyle
\parallel \widetilde q \parallel \ib{0}^{2}\quad \geq \quad
{1\over 6} \,\,  \sum_{j=0}^{n-1} h_{j+1/2} \,\, (q_{j}^{2}+q_{j+1}^{2})\,. $

\smallskip \noindent 
In an analogous way, we have 

\smallskip \noindent  $ \displaystyle
\parallel q \parallel \ib{0}^{2} \,\, \leq \,\,  \sum_{j=0}^{n-1} h_{j+1/2} \,
(\mid q_{j} \mid + \mid q_{j+1}  \mid ) ^{2} \,\int_{0}^{1}
(\psi(\theta))^{2}\,{\rm d}\theta $

\smallskip \noindent 
{ \it i.e.}

\smallskip \noindent  (4.20) $\qquad \displaystyle
\parallel q \parallel \ib{0}^{2} \quad \leq \quad 2 \,\,
\widetilde{\delta}\,\,\,  \sum_{j=0}^{n-1} h_{j+1/2}  \,\,
 \bigl(\mid q_{j} \mid  ^{2}+ \mid q_{j+1}  \mid  ^{2} \bigr) \,. \,$
 
\smallskip \noindent
We have the same inequality when the interpolant function $q$ is replaced by
$\widetilde q$, and $\delta$ replaced by its value when $\psi(\smb)$ is
replaced by affine interpolation $\theta \mapsto \theta$ :

\smallskip \noindent  (4.21) $\qquad \displaystyle
\parallel \widetilde q \parallel \ib{0}^{2} \quad\leq\quad 
{2\over3} \,\, \sum_{j=0}^{n-1} h_{j+1/2}  \,\,
 (\mid q_{j} \mid  ^{2}+ \mid q_{j+1}  \mid  ^{2})  \,. $

\smallskip \noindent $\bullet \quad$
From (4.20) and (4.19) we deduce

\smallskip \noindent  $ \displaystyle
\parallel q \parallel \ib{0}^{2} \quad\leq\quad 2 \,\,
\widetilde{\delta}\,\, \sum_{j=0}^{n-1} h_{j+1/2}  \,\,
 (\mid q_{j} \mid  ^{2}+ \mid q_{j+1}  \mid  ^{2}) \quad\leq\quad
 12\,\, \widetilde{\delta}\,\, \parallel \widetilde q \parallel \ib{0}^{2} \,$
 
\smallskip \noindent
that establishes the second inequality of (4.14). Using estimates (4.18) and
(4.21) we have

\smallskip \noindent  $ \displaystyle
\parallel q \parallel \ib{0}^{2} \quad\geq\quad 
\delta\,\, \sum_{j=0}^{n-1} h_{j+1/2}  \,\,
 (\mid q_{j} \mid  ^{2}+ \mid q_{j+1}  \mid  ^{2}) \quad\geq\quad
 {3\over2}\,\, \delta\,\, \parallel \widetilde q \parallel \ib{0}^{2} \,$
 
\smallskip \noindent
and the proof of inequality (4.14) is completed. $\hfill  \square \kern0.1mm$

\bigskip \noindent $\bullet \quad$
We show now that if condition (4.7) of compatibility interpolation condition is
not satisfied, the uniform inf-sup condition (4.6) cannot be satisfied for any
family of uniformly regular meshes. In other words, trial functions in space 
$\,\QT^{\psi}\,$ oscillate too
much and stability is in defect.  

\bigskip
\noindent {\bf Theorem 3. $\quad$  Lack of inf-sup condition.}

\noindent 
Let $\psi : [0,1] \rightarrow \R\quad$ be a continuous function satisfying
conditions  (3.4), (3.7) and the negation of compatibility interpolation
condition, {\it i.e.} 

\smallskip \noindent  (4.22) $\qquad \displaystyle
\exists \,\theta \in\, ]0,1[\,,\quad\psi(\theta)+\psi(1\!-\!\theta)\,\ne 1\,.$

\smallskip \noindent  
Then for any family $\cal{U}_{\alpha,\beta}$ of uniformly regular meshes
($0<\alpha<1<\beta$), the inf-sup condition (4.6) is not satisfied for spaces
$V_{1}=\UT\times\PT$ and $V_{2}=\UT\times\QT^{\psi}$ and  meshes $\cal{T}$ of 
$\cal{U}_{\alpha,\beta}$:

\setbox21=\hbox{$\displaystyle \forall (\alpha,\beta)\,,\,\,
0<\alpha<1<\beta\,,\quad  \forall\,D>0\,,$}
\setbox22=\hbox{$\displaystyle  \exists\,\cal{T} \in
\cal{U}_{\alpha,\beta}\,,\quad \exists\,\xi\in \UT \times \PT  \,$ such that $
\,\,\,\,\parallel \xi \parallel \,\,=\,\, 1\,\, $ and }
\setbox23=\hbox{$ \displaystyle  \forall \,\eta \in
\UT\times\QT^{\psi}\,,\quad \parallel
\eta \parallel \,\leq\,  1 \quad \Rightarrow \quad 
\, \gamma(\xi,\eta) \, \leq \, D\,.$}
\setbox30= \vbox {\halign{#\cr \box21 \cr    \box22 \cr   \box23 \cr  }}
\setbox31= \hbox{ $\vcenter {\box30} $}
\setbox44=\hbox{\noindent  (4.23) $\displaystyle  \quad  \left\{ \box31 \right.
$}  

\smallskip \noindent $ \box44 $

\bigskip \noindent  {\bf Proof of theorem 3.}

\smallskip \noindent $\bullet\quad$ 
The first point what we have to show is that if relation (4.22) is satisfied,
then we have
\setbox21=\hbox{$\displaystyle \int_{0}^{1}
{{{\rm d}\psi}\over{{\rm d}\theta}} (\theta)\,{{{\rm d}\psi}\over{{\rm d}\theta}}
(1\!-\!\theta)\,{\rm d}\theta$}

\smallskip \noindent  (4.24) $\qquad \displaystyle
\mod{\box21} \quad < \quad \int_{0}^{1} \Bigl( {{{\rm d}\psi}\over{{\rm d}\theta}}
\Bigr)^{2} \, {\rm d}\theta \,. $

\smallskip \noindent 
The large inequality between the two sides of (4.24) just express Cauchy-Schwarz
inequality. If the equality is realized,  functions
$\displaystyle \,{{{\rm d}\psi}\over{{\rm d}\theta}}(\smb)\,$ and
$\displaystyle \,{{{\rm d}\psi}\over{{\rm d}\theta}}(1\!-\!\smb)\,$ are 
linearly dependent:

\smallskip \noindent  (4.25) $\quad \displaystyle
\exists\,(\lambda,\mu)\in \R\,,\,\, (\lambda,\mu)\ne(0,0)\,,
\quad\forall\,\theta\in[0,1]\,,\quad \lambda \,  {{{\rm
d}\psi}\over{{\rm d}\theta}}(\theta) - \mu\, {{{\rm d}\psi} \over{{\rm
d}\theta}}(1\!-\!\theta)\,=\,0\,. \,$

\smallskip \noindent 
Then function $\,\,\theta \mapsto \lambda \, \psi(\theta) +\mu\,
\psi(1\!-\!\theta)\,\,$is equal to some constant whose value is equal to
$\mu$ (take $\theta=0$ and apply (3.4)). Moreover, taking $\theta=1$, we get
$\,\lambda=\mu\,$ and we obtain in this way

\smallskip \noindent  (4.26) $\qquad \displaystyle
\mu\,\bigr( \psi(\theta) + \psi(1\!-\!\theta) - 1 \bigr)\,\,=\,\,
0\,,\quad \forall\,\theta \in [0,1] \,.$

\smallskip \noindent 
Joined with relation (4.22), $\mu$ is necessarily equal to zero and finally
$\,\lambda=\mu\,=\,0$ which express the contradiction.

\monitem 
We set

\smallskip \noindent  (4.27) $\qquad \displaystyle
\epsilon \quad=\quad \int_{0}^{1} \Bigl( {{{\rm d}\psi}\over{{\rm d}\theta}}
\Bigr) ^{2}\,{\rm d}\theta \,-\, \int_{0}^{1}  {{{\rm d}\psi}
\over{{\rm d}\theta}}(\theta)\, {{{\rm d}\psi}\over{{\rm d}\theta}} (1
\!-\!\theta) \,{\rm d}\theta \,$

\smallskip \noindent 
and $\,\epsilon > 0\,$ due to (4.24). We evaluate now the $L^{2}$ norm of
$\displaystyle \,{\rm div}\,q\,=\,{{{\rm d}q}\over{{\rm d}x}}\,=\,{{{\rm
d}}\over{{\rm d}x}}\,\bigl( \sum_{j=0}^{n-1}q_{j}\,\psi_{j}\,\bigr)\,$:

\setbox20=\hbox{$\parallel {\rm div}\,q \parallel \ib{0}^{2}$}
\setbox21=\hbox{$\displaystyle  \sum_{j=0}^{n-1} \, \int_{
\displaystyle x_{j}}^{\displaystyle x_{j+1}}
\biggl(  \, {{\rm d}\over{{\rm d}x}}\,\Bigl[\, 
q_{j} \,\, \psi\Bigl({{x_{j+1}-x}\over{h_{j+1/2}}}\Bigr)
+ q_{j+1} \,\, \psi \Bigl({{x-x_{j}}\over{h_{j+1/2}}}\Bigr) \Bigr]
\,\biggr)^{2} \,{\rm d}x$}
\setbox22=\hbox{$\displaystyle  \sum_{j=0}^{n-1} {{1}\over{h_{j+1/2}}} 
\int_{0}^{1} \Bigl(  -q_{j} \, {{{\rm d}\psi}\over{{\rm d}\theta}}(1\!-\!\theta)
\,+\,  q_{j+1} \, {{{\rm d}\psi}\over{{\rm d}\theta}}(\theta) \Bigr)^{2} \,
{\rm d}\theta $}
\setbox23=\hbox{$\displaystyle  \sum_{j=0}^{n-1}  {{1}\over{h_{j+1/2}}}  
\biggl[ (q_{j}^{2}+q_{j+1}^{2})\,\int_{0}^{1} \Bigl({{{\rm d}\psi} \over {{\rm
d}\theta}}(\theta)\Bigr)^{2}\,{\rm d}\theta \biggr]$}
\setbox24=\hbox{$\displaystyle  \qquad \qquad -2 \,\,\sum_{j=0}^{n-1} 
{{1}\over{h_{j+1/2}}}   \biggl[  \,q_{j}\,q_{j+1} \int_{0}^{1}
{{{\rm d}\psi}\over{{\rm d}\theta}}(\theta)\, {{{\rm d}\psi}\over{{\rm d}\theta}}
(1\!-\!\theta) \,{\rm d}\theta \bigg]$}
\setbox25=\hbox{$\displaystyle \int_{0}^{1} {{{\rm d}\psi} \over {{\rm
d}\theta}}(\theta)\, {{{\rm d}\psi}\over{{\rm d}\theta}}(1\!-\!\theta)
\,{\rm d}\theta$}
\setbox26=\hbox{$\displaystyle \epsilon\, \sum_{j=0}^{n-1}
{{(q_{j}^{2}+q_{j+1}^{2})}\over{h_{j+1/2}}} + \sum_{j=0}^{n-1}
{{(\mid q_{j}\mid \! -\!\mid q_{j+1} \mid )^{2}}\over{h_{j+1/2}}}
\,\mod{\box25}\,.$}
\setbox27=\vbox{\halign{#&#&#\cr \box20&$\quad=\quad$&\box21 \cr 
&$\quad=\quad$& \box22 \cr &$\quad=\quad$& \box23 \cr 
&$\quad \, \quad$& \box24  \cr &$ \quad \geq \quad $& \box26  \cr}}
\smallskip \noindent  $\box27$

\smallskip \noindent
Then

\smallskip \noindent  (4.28) $\qquad \displaystyle
\parallel {\rm div}\,q \parallel \ib{0}^{2} \quad \geq \quad \epsilon\,
\sum_{j=0}^{n-1} {{(q_{j}^{2}+q_{j+1}^{2})}\over{h_{j+1/2}}} \,.$

\monitem 
We establish now (4.23) which express the negation of uniform inf-sup condition.
Consider a mesh $\cal T$ composed with $n$ elements uniformly distributed :

\smallskip \noindent  $ \displaystyle
0 \,\,=\,\, x_{0} < x_{1}={{1}\over{n}} <\cdots<x_{k} \,\,=\,\,
{{k}\over{n}}<\cdots<x_{n-1}<x_{n} \,\,=\,\, 1 \,$

\smallskip \noindent
with integer $n$ chosen such that

\smallskip \noindent  (4.29) $\qquad \displaystyle
{2\over{\sqrt{n\epsilon}}}\quad\leq\quad D\,. \,$

\smallskip \noindent
It is clear that for each pair $(\alpha,\beta)$ satisfying relation (4.10), mesh
$\cal T$ defined previously belongs to  $\cal{U}_{\alpha,\beta}$
($h_{j+1/2}^{\cal{T}}$ is exactly equal to $\,\displaystyle {{1}\over{\nT}}\,$
with notations proposed at Definition 1). Introduce $\,u(x)\equiv\,1$,
$\,p(x)\equiv\,0$ and $\xi\equiv(u,p)=(1,0)\,$ which is clearly of norm
equal to unity in space $\,V=L^{2}(0,1)\times H^{1}(0,1)\,$. For each
$\,\eta=(v,q)\,$ in subspace $\UT\times\QT^{\psi}$, we have

\smallskip \noindent  $ \displaystyle
\gamma(\xi,\eta) \quad =\quad (1\,,\,{\rm div}\,q) \quad =\quad
q(x_{n})\,-\,q(x_{0})\,. \,$

\smallskip \noindent
From inequality (4.28) we have :

\smallskip \noindent  $ \displaystyle
\mid q_{j} \mid^{2}\quad\leq\quad {{1}\over{\epsilon}}\,h_{j+1/2}\, \parallel
{\rm div}\,q \parallel \ib{0}^{2}\quad\leq\quad
{1\over{n\,\epsilon}}\,,\quad\forall\,j=0,\cdots,n \,$

\smallskip \noindent
when $\cal T$ is chosen as above and $\eta$ with a norm less or equal to $1$ in
space $L^{2}(0,1)\times H^{1}(0,1)\,$ (see (1.29)). Then we have

\smallskip \noindent  (4.30) $\qquad \displaystyle
\mid \gamma(\xi,\eta) \mid \quad \leq\quad
{2\over{\sqrt{n\epsilon}}}\quad\leq\quad D \,$

\smallskip \noindent
if relation (4.29) is realized. Relation (4.23) is proven and uniform inf-sup
condition is in defect.   $\hfill  \square \kern0.1mm$

\bigskip \bigskip \newpage   \noindent {\bf \large 5) \quad  Convergence of finite volumes  in
the one dimensional case} 

\monitem  
We have proven in section 4 (Theorem 3) that if the compatibily interpolation
condition

\smallskip \noindent  (5.1) $\qquad \displaystyle
\psi(\theta) + \psi(1-\theta) \,\equiv\, 1 \,,\qquad \forall \,\theta \,
\in [0,1] \,$

\smallskip \noindent
is not realized, there is no hope to obtain convergence in usual Hilbert spaces
for the finite volume method  (1.5)-(1.9) formulated  as a mixed Petrov-Galerkin
finite element method (3.1)-(3.3) associated with a family 
$\cal{U}_{\alpha,\beta}$ of uniformly reguler meshes $\cal T$, shape functions
$\,\xiT\!\!=(\uT,\pT) \in \UT \times \PT \,$, weighting functions
$ \eta = (v,q)\in \UT \times \QT^{\psi}$ and bilinear form

\smallskip \noindent  (5.2) $\qquad \displaystyle
\gamma(\xiT,\eta)\quad=\quad(\pT,q) \,+\, (\uT,{\rm div}\,q)\,+\,
({\rm div}\,\pT,v)\,.\,$

\smallskip \noindent 
On the contrary, if compatibility interpolation condition (5.1) is realized, we
have convergence and the following result holds.

\bigskip
\noindent {\bf Theorem 4. $\quad$  Convergence of 1D finite volumes in Hilbert
spaces.}

\noindent 
Let  $\psi : [0,1] \rightarrow \R\quad$ be a continuous function, satisfying
$\,\psi(0)=0\,$, the compatibility interpolation condition (5.1) and
orthogonality condition

\smallskip \noindent  (5.3) $\qquad \displaystyle
\int_{0}^{1} (1\!-\!x)\,\psi(x)\,{\rm d}x\quad=\quad0\,.\,$

\smallskip \noindent  
Let $\cal{U}_{\alpha,\beta}$ $(0<\alpha<1<\beta\,)$ be a family af regular meshes
$\cal T$ in the sense given in definition 2, $\UT$ and $\PT$ be interpolation
spaces of piecewise constant functions in each element and continuous
piecewise linear functions, $\QT^{\psi}$ be the space of weighting functions
proposed at Definition 1 : function $\psi_{j}$ is defined in (3.5) and
function $\,q\in  \QT^{\psi}$ satisfies 

\smallskip \noindent  (5.4) $\qquad \displaystyle
q\quad=\quad \sum_{j=0}^{n} q_{j}\, \psi_{j}\,.\,$

\smallskip \noindent  
Then for each $\,f\in L^{2}\,$, the solution $\xiT=(\uT,\pT)\in \UT\times \PT$ of
the finite volume method for the approximation of the solution
$\,\xi\equiv(u,p={\rm grad}\,u)\,$ of Dirichlet problem for one-dimensional 
Poisson equation

\smallskip \noindent   (5.5) $\qquad \displaystyle
-\Delta\,u\,=\,f\,\,\,{\rm in\,\,\,}]0,1[\,,\qquad u(0)\,=\,
u(1)\,=\,0\, $

\smallskip \noindent 
is given by solving problem (3.1)-(3.3) :

\smallskip \noindent   (5.6) $\qquad \displaystyle
\gamma (\xiT,\eta)\,=\, (f,v) \,,\qquad \forall\,\eta\,=\, (v,q)\in
\UT\times \QT^{\psi} \,$

\smallskip \noindent
where bilinear form $\gamma(\smb,\smb)$ is defined in (5.2).

\noindent 
Moreover when $f$ belongs to space $H^{1}(0,1)$, there exists some constant $C>0$
depending only on $\alpha$ and $\beta$ such that

\smallskip \noindent   (5.7) $\qquad \displaystyle
\parallel u-\uT \parallel \ib{0} + \parallel p-\pT \parallel \ib{1} \quad
\leq\quad C\,\,\hT \,\parallel f \parallel \ib{1}\,,\qquad \forall \,\cal{T} \in
\cal{U}_{\alpha,\beta} \,$

\smallskip \noindent
where $\hT$ is the maximal step size of mesh $\cal{T}$ precisely defined in
(4.13).

\bigskip
\noindent {\bf Remark 1.}
A simple but fundamental remark is that the finite volume method (1.5)-(1.9)
corresponds exactly to the mixed Petrov-Galerkin finite element formulation,
independently of the choice of interpolation function $\psi$ satisfying (5.1).
This is due to the fact that the heuristic relation (1.5) holds if the following
relation

\smallskip \noindent   (5.8) $\qquad \displaystyle
\int_{0}^{1}\theta\,\psi(\theta)\, {\rm d}\theta \quad=\quad{1\over2} \,$

\smallskip \noindent 
(see also (3.9) and  (3.11)) is satisfied. But relation (5.8)  derives clearly
from relation (5.3) by integration of identity (5.1) after multiplication by
$\theta$. Compatibily interpolation condition (5.1) gives an acute link between
consistency (relation (5.8)) and convergence (inf-sup condition (4.6)). We have
proven that the heuristic relation (1.5) is the only  possible finite volume
scheme associated with a stable mixed Petrov Galerkin  formulation. 

\bigskip \noindent $\bullet \quad$
Some propositions are usefull to be established, before prooving completely
Theorem~4, first established with other techniques by Baranger {\it et al} [BMO96] 
and also studied with finite difference techniques by Eymard,  Gallou\"et and
Herbin [EGH2k].

\bigskip
\noindent {\bf Proposition 5. $\quad H^{1}$ continuity of $P_{1}$
interpolation.}

\noindent 
Let $\PiT$ be the classical  $P_{1}$ interpolation operator in space $\PT$,
defined by

\smallskip \noindent   (5.9) $\qquad \displaystyle
\bigl( \PiT\mu\bigr)(x_{j}) \,\,=\,\,\mu(x_{j})\,,\quad \forall \,\mu \in
H^{1}(0,1)\,,\quad \forall\,x_{j}\,{\rm  vertex \,\,of\,\, mesh
\,\,}\cal{T}\,.$

\smallskip \noindent 
When mesh $\,\cal{T}\,$ describes  a family  $\,\cal{U}_{\alpha,\beta}\,$ of
uniformly regular meshes, we have the following property :

\smallskip \noindent   (5.10) $\qquad \displaystyle
\exists\,C_{1}>0\,,\quad\forall\,\mu\in H^{1}(0,1)\,,\quad \parallel 
\PiT\mu \parallel \ib{1} \,\leq\, C_{1}\, \parallel \mu \parallel \ib{1}\,.$

\bigskip \noindent {\bf Proposition 6. $\quad$ Discrete stability.}

\noindent 
Let $\alpha$ and $\beta$ be such that $0<\alpha<1<\beta$ and
$\,\cal{U}_{\alpha,\beta}\,$ be a family of uniformly regular meshes. When
$\,\psi\,$ is chosen satisfying hypotheses of Theorem~4, there exists some
constant
$\,C>0\,$ such that

\setbox21=\hbox{$\quad \forall\,{\cal{T}}\!\! \in {\cal{U}}_{\alpha,\beta} \,
\quad\forall\, u \in \UT\, \quad  \exists\,q\in \QT^{\psi},\,\,$}
\setbox22=\hbox{\quad $(u \,,\,{\rm div}\,q)\,\,=\,\,\parallel u \parallel
\ib{0}^{2}\quad{\rm and\quad}
\parallel q \parallel \ib{1}\,\,\leq\,\,C\,\parallel u \parallel \ib{0}\,.$}
\setbox30= \vbox {\halign{#\cr \box21 \cr   \box22 \cr   }}
\setbox31= \hbox{ $\vcenter {\box30} $}
\setbox44=\hbox{\noindent  (5.11) $\displaystyle  \quad   \left\{ \box31 \right.
$}  

\smallskip \noindent $ \box44 $

\bigskip \noindent  {\bf Proof of proposition 6.}
\smallskip \noindent $\bullet \quad$
Let $u$ be given in $\UT$ and $\,\varphi\in H_{0}^{1}(0,1)\,$ be the variational
solution of the problem

\smallskip \noindent   (5.12) $\qquad \displaystyle
\Delta \chi\,\,=\,\,u\,\,\,{\rm on}\,\,
]0,1[\,,\quad \quad \quad \chi(0)\,=\chi(1)\,=\,0\,.$

\smallskip \noindent
Then (see {\it e.g.} [Ad75]), $\chi$ belongs to space $H^{2}$ and there exists
some constant $C_{2}$ independent on $u$ such that

\smallskip \noindent   $ \displaystyle
\parallel \chi \parallel \ib{2}\,\,\leq\,\, C_{2}\,\parallel u
\parallel \ib{0}\,. $

\smallskip \noindent 
Let $\,\widetilde{q}\,=\,\PiT({\rm grad}\,\chi)\,$ be the usual $P_{1}$
interpolate of ${\rm grad}\,\chi$. From Proposition 5, we have

\smallskip \noindent   (5.13) $\qquad \displaystyle
\parallel \widetilde{q} \parallel \ib{1}\,\leq\,C_{1}\,\parallel {\rm grad}\,\chi
\parallel \ib{1}\,\leq\,C_{1}\,\parallel \chi \parallel \ib{2}
\,\leq\, C_{1}\,C_{2} \parallel u \parallel \ib{0} \,=\, 
C_{3} \parallel u \parallel \ib{0} \,. $

\smallskip \noindent  
Writing $\,\widetilde q ={\sum_{j=0}^{n}} q_{j} \,\varphi_{j} \in \PT\,$, we
introduce the second interpolant function $\,  q ={\sum_{j=0}^{n}} q_{j}
\,\psi_{j} \,\, \in \QT^{\psi}\,\,\,$ and we have, for any $\,v \in\UT$

\setbox20=\hbox{$({\rm div}\,q \,,\, v)$}
\setbox21=\hbox{$\displaystyle  \sum_{j=0}^{n-1} v_{j+1/2} \,
\int_{\displaystyle x_{j}}^{\displaystyle x_{j+1}} {\rm div}\,q\,\,{\rm d}x$}
\setbox22=\hbox{$\displaystyle \sum_{j=0}^{n-1} v_{j+1/2} \,\bigl( q_{j+1}-q_{j}
\bigr) $}
\setbox23=\hbox{$\displaystyle   \sum_{j=0}^{n-1} v_{j+1/2} \,\Bigl(
{{{\rm d}\chi}\over{{\rm d}x}}(x_{j+1}) -{{{\rm d}\chi}\over{{\rm d}x}}(x_{j}) \Bigr) $}
\setbox24=\hbox{$\displaystyle   \sum_{j=0}^{n-1} v_{j+1/2}
\,\int_{\displaystyle x_{j}}^{\displaystyle x_{j+1}} \Delta \chi \,{\rm d}x$}
\setbox25=\hbox{$\displaystyle  (u\,,\, v)\,,\qquad\forall\,v\in \UT\,. $}
\setbox27=\vbox{\halign{#&#&#\cr \box20&$\quad=\quad$&\box21 \cr 
&$\quad=\quad$& \box22 \cr &$\quad=\quad$& \box23 \cr  &$\quad=\quad$& \box24 \cr
 &$\quad=\quad$& \box25 \cr}}

\smallskip \noindent  $\box27$ 

\smallskip \noindent 
In particular (choose $\,v=u\,)$, the equality $\,(u \,,\, {\rm
div}\,q)=\,\parallel u \parallel \ib{0}^{2}\,$ of relation (5.11) is established.

\bigskip \noindent $\bullet \quad$
We show now the stability inequality of relation (5.11), between  $\,\parallel q
\parallel \ib{1}\,$ and $\,\parallel u \parallel \ib{0}\,$. We have, from
relation (4.14) of Proposition 4 and estimations (5.13)

\smallskip \noindent   $ \displaystyle
\parallel q \parallel \ib{0}^{2} \quad\leq\quad12 \, \widetilde{\delta} \,
\parallel \widetilde q \parallel \ib{0}^{2}  \quad\leq\quad 12 \,\,
\widetilde{\delta}\,\, C_{3}^{2}\, \parallel u \parallel \ib{0}^{2} \,$

\smallskip \noindent
and since (5.1) holds,
\setbox20=\hbox{$\mid q \mid \ib{1}^{2}$}
\setbox21=\hbox{$\displaystyle  \sum_{j=0}^{n-1} {{1}\over{h_{j+1/2}}} \,
(q_{j+1}-q_{j})^{2}\,\int_{0}^{1} \Bigl( {{{\rm d}\psi}\over{{\rm d}\theta}}
\Bigr)^{2}\,{\rm d}\theta$}
\setbox22=\hbox{$\displaystyle \int_{0}^{1} \Bigl( {{{\rm d}\psi}\over{{\rm d}\theta}}
\Bigr)^{2}\,{\rm d}\theta \quad \mid \widetilde{q} \mid \ib{1}^{2} $}
\setbox23=\hbox{$\displaystyle  \int_{0}^{1} \Bigl( {{{\rm d}\psi}\over{{\rm d}\theta}}
\Bigr)^{2}\,{\rm d}\theta \quad  C_{3}^{2}\,\, \parallel u \parallel \ib{0}^{2}\,.$}
\setbox27=\vbox{\halign{#&#&#\cr \box20&$\quad=\quad$&\box21 \cr 
&$\quad=\quad$& \box22 \cr &$\quad \leq \quad$& \box23 \cr }}
\smallskip \noindent  $\box27$ 

\smallskip \noindent
From these inequalities, we deduce inequality $\, \parallel q \parallel
\ib{1}\,\,\leq\,\,C\,\parallel u \parallel \ib{0}\,$, with

\smallskip \noindent   $ \displaystyle
C\quad=\quad \biggl(12 \, \widetilde{\delta}\,\,+\,\int_{0}^{1} \Bigl(
{{{\rm d}\psi}\over{{\rm d}\theta}} \Bigr)^{2}\,{\rm d}\theta\,
\biggr)^{1/2}\,\,C_{3}\, $

\smallskip \noindent
and Proposition 6 is established.  $\hfill  \square \kern0.1mm$

\bigskip
\noindent {\bf Proposition 7. $\quad$ Uniform discrete inf-sup condition.}

\noindent
Let $\psi : [0,1] \rightarrow \R\quad$ be a continuous function satisfying
$\,\psi(0)=0\,$, orthogonality condition (5.3) and compatibility interpolation
condition (5.1). Let $\widetilde{\delta}$ be defined according to relation (4.16)
and

\smallskip \noindent   (5.14) $\qquad \displaystyle
K\,\,=\,\,{4\over3}\,\biggl( 1+\sqrt{12\, \widetilde{\delta}}\,
\biggr)\,.$

\smallskip \noindent
Let $\,\alpha$ and $\beta\,$ be real numbers such that $\,0<\alpha<1<\beta\,$, 
$\cal{U}_{\alpha,\beta}$ be a family of uniformly regular meshes, 
$\,\gamma(\smb,\smb)\,$ be the bilinear form defined in (5.2), $\,C\,$
be the constant associated with inequality (5.11) in Proposition 6 and
$\,\rho>0\,$ be chosen such that

\smallskip \noindent   (5.15) $\qquad \displaystyle
\rho\,+\sqrt{K\,\rho}\quad\leq\quad{1\over{C}}\,\,\sqrt{1-\rho^{2}-K\rho}\,. $

\smallskip \noindent  
Then we have the following uniform discrete inf-sup condition :


\setbox21=\hbox{$\quad \forall\,{\cal{T}}\!\!\in
{\cal{U}}_{\alpha,\beta}\,,\quad\forall\, \xi \, =\, (u,p) \in
\UT\times\PT\,,  \parallel \xi \parallel \, =1\,,  $} 

\setbox22=\hbox{$\quad \exists\,\eta=(v,q)\in\UT\times\QT^{\psi}\,,\quad 
\parallel 
\eta \parallel\,\leq\,1 \quad {\rm and}\quad
\gamma(\xi,\eta)\,\,\geq\,\,\rho\,.$}


\smallskip \noindent   (5.16) $\qquad \displaystyle 
\left\{ \begin{array} {l}  \displaystyle  
 \box21 \\ \box22 
 \end{array} \right.  $ 

\bigskip \noindent  {\bf Proof of proposition 7.}

\monitem 
As in Proposition 1, we distinguish between three cases. If we have the condition

\smallskip \noindent   (5.17) $\qquad \displaystyle
\parallel {\rm div}\,p \parallel \ib{0} \quad\geq\quad \rho\,,$

\smallskip \noindent  
let $\,\eta\equiv(v,q)\,$ be defined by
$\displaystyle \,\, v={{{\rm div}\,p}\over{\parallel {\rm div}\,p
\parallel \ib{0}}}\,$ and $\,\,q=0\,$. Then, due to relation (5.2), we have
$\gamma(\xi,\eta) = ({\rm div}\,p,v)=\,\parallel {\rm div}\,p \parallel \ib{0}\,\geq\,\rho\,$
and inequality (5.16) is proven in this simple case.

\smallskip  \noindent $\bullet \quad$
When (5.17) is in defect, we suppose also that $\,p\,$ is sufficiently large :

\smallskip \noindent   (5.18) $\qquad \displaystyle
\parallel {\rm div}\,p \parallel \ib{0} \quad\leq\quad \rho\quad{\rm
and}\quad \parallel p \parallel \ib{0}^{2} \quad\geq\quad K\,\rho \,. $

\smallskip  \noindent
We set $\,p=\sum_{j=0}^{n}\,p_{j} \varphi_{j}\,$ and introduce $\,q\in
\QT^{\psi}\,$ according to the relation

\smallskip \noindent   (5.19) $\qquad \displaystyle
q\quad=\quad{1\over{\sqrt{12 \, \widetilde{\delta}}}}\, \sum_{j=0}^{n} p_{j}
\psi_{j}\,. $ 

\smallskip \noindent 
From inequality (4.14) and the hypothesis done on $\,\xi=(u,p)\,$, we have

\smallskip \noindent  $ \displaystyle
\,\parallel q \parallel\ib{0} \quad \leq\quad \parallel p \parallel
\ib{0}\quad \leq \quad \, 1 $

\smallskip \noindent
and moreover :

\setbox20=\hbox{$(p,q)$}
\setbox21=\hbox{$\displaystyle  {1\over{\sqrt{12 \, \widetilde{\delta}}}}\,
\sum_{j=0}^{n} p_{j}^{2}\,(\varphi_{j},\psi_{j})$}
\setbox22=\hbox{$\displaystyle {1\over{\sqrt{12 \, \widetilde{\delta}}}}\,
\sum_{j=0}^{n-1}  h_{j+1/2}\,\,\bigl( p_{j}^{2} + p_{j+1}^{2} \bigr)\,\,
\int_{0}^{1} \theta \, \psi(\theta)\,{\rm d}\theta \quad $ due to (3.9) and
(1.8) }
\setbox23=\hbox{$\displaystyle  {1\over{4 \sqrt{3\, \widetilde{\delta}}}}\,\,
\sum_{j=0}^{n-1} h_{j+1/2}\,\,\bigl( p_{j}^{2} + p_{j+1}^{2} \bigr)$}
\setbox24=\hbox{$\displaystyle 
{1\over8}\,\sqrt{{{3}\over{\widetilde\delta}}}\,\,
\parallel p \parallel \ib{0}^{2} \qquad {\rm due \,\,to \,\,(4.21).}$}
\setbox27=\vbox{\halign{#&#&#\cr \box20&$\quad=\quad$&\box21 \cr 
&$\quad=\quad$& \box22 \cr &$\quad = \quad$& \box23 \cr 
&$\quad\geq\quad$& \box24 \cr}}

\smallskip \noindent $\box27$

\smallskip \noindent
We introduce  $\,\eta=(0,q)\,$. Then we have shown that $\,\parallel \eta
\parallel \, \leq\,1\,$ and we have also
\setbox20=\hbox{$\gamma(\xi,\eta)$}
\setbox21=\hbox{$\displaystyle (p,q) + (u,{\rm div}\,q)$}
\setbox22=\hbox{$\displaystyle (p,q)  + 
{1\over{\sqrt{12 \, \widetilde{\delta}}}} \,(u,{\rm div}\,p)$}
\setbox23=\hbox{$\displaystyle  {1\over8}\,\sqrt{{{3}\over{\widetilde\delta}}}
\, K \rho \,\, - \,\,{1\over{\sqrt{12 \, \widetilde{\delta}}}} \,\rho$}
\setbox24=\hbox{$\displaystyle \rho$}
\setbox27=\vbox{\halign{#&#&#&#&# \cr \box20&$\quad=\quad$&\box21&& \cr 
&$\quad=\quad$& \box22&& \cr 
&$\quad \geq \quad$& \box23&$\quad=\quad$&\box24\cr}}

\smallskip \noindent $\box27$

\smallskip \noindent
due to (5.14). Then (5.16) holds in this second case. 

\smallskip \noindent $\bullet \quad$
In the third case, we suppose

\smallskip \noindent   (5.20) $\qquad \displaystyle
\parallel {\rm div}\,p \parallel \ib{0} \, \, \leq \,\, \rho\,,\quad \parallel p
\parallel \ib{0}^{2}  \, \, \leq \,\,  K\,\rho\,. $

\smallskip \noindent  
Then because the norm of $\xi$ is exactly equal to $1$, we have

\smallskip \noindent   $ \displaystyle
\parallel u \parallel \ib{0}^{2}\,\,=\,\, 1 \,- \parallel p \parallel \ib{0}^{2}
- \parallel {\rm div} \,p \parallel \ib{0}^{2}\quad\geq\,\, 1-K\,\rho-\rho^{2} $

\smallskip \noindent 
which is strictly positive because the right hand side of inequality (5.15) is
strictly positive ($\rho > 0$). Let $\,q\,$ be associated with $u$ according to
relation (5.11) of proposition 6 :

\smallskip \noindent   (5.21) $\qquad \displaystyle
q\in \QT^{\psi}\,,\quad (u \,,\, {\rm div}\,q)\,\,=\,\,\parallel u \parallel
\ib{0}^{2}\,,\quad \parallel q \parallel \ib{1}\,\,\leq\,\,C\,\parallel u
\parallel \ib{0}\,. $

\smallskip \noindent 
Then $ \displaystyle \,\eta\equiv(0,{1\over{C\,\parallel u \parallel
\ib{0}}}\,q)\,$ has a norm not greater than $1$ and due to relation (5.2), we
have 
\setbox20=\hbox{$\gamma(\xi,\eta)$}
\setbox21=\hbox{$\displaystyle \Bigl(  p \,,\,{{q}\over{C\,\parallel u \parallel
\ib{0}}} \Bigr)  + \Bigl( u\,,\,{{{\rm div}\,q}\over{C\,\parallel u \parallel
\ib{0}}} \Bigr) $}
\setbox22=\hbox{$\displaystyle -\, \parallel p \parallel \ib{0}\,\,+\,\,
{1\over{C}}\,\parallel u \parallel \ib{0}$}
\setbox42=\hbox{$\qquad$ due to (5.21)}
\setbox23=\hbox{$\displaystyle -\,\sqrt{K\,\rho}\,\,+\,\, 
{1\over{C}}\,\,\sqrt{1-\rho^{2}-K\rho}$}
\setbox43=\hbox{$\qquad$ due to (5.20)}
\setbox24=\hbox{$\displaystyle \rho$}
\setbox44=\hbox{$\qquad$ due to (5.15)}
\setbox27=\vbox{\halign{#&#&#&# \cr \box20&$\quad=\quad$&\box21 \cr 
&$\quad\geq\quad$& \box22 & \box42  \cr &$\quad \geq \quad$& \box23 & \box43 \cr 
&$\quad\geq\quad$& \box24 & \box44 \cr}}

\smallskip \noindent $\box27$

\smallskip \noindent 
that ends the establishment of uniform inf-sup condition  (5.16).
 $\hfill  \square \kern0.1mm$
 
\bigskip \noindent $\bullet \quad$
We need also interpolation results, that are classical (see, {\it e.g.} [CR72]).
We detail them for completeness.

\bigskip
\noindent {\bf Proposition 8. $\quad$ Interpolation errors.}

 \noindent 
Let $\,v\in L^{2}(0,1)\,$ and $\,q \in H^{1}(0,1)\,$ be two given functions, 
$\,\MT\,$ and $\,\PiT\,$ the piecewise constant ($P_{0}$) and continuous
piecewise linear  ($P_{1}$) interpolation operators on mesh $\,\,\cal{T}\,$
defined in finite dimensional spaces $\,\UT\,$ and $\,\PT\,$
respectively by the following relations

\smallskip \noindent   (5.22) $\qquad \displaystyle
\bigl( \MT v\bigr)(x) \,\,=\,\,{{1}\over{h_{j+1/2}}}\,\int_{\displaystyle
x_{j}}^{\displaystyle x_{j+1}} v(y)\,{\rm d} y\,,\quad x_{j}<x<x_{j+1}  \,$

\smallskip \noindent   (5.23) $\qquad \displaystyle
\bigl( \PiT q \bigr)(x) \,\,=\,\,q(x_{j})\,{{x_{j+1}-x}\over{h_{j+1/2}}} + 
q(x_{j+1})\,{{x-x_{j}}\over{h_{j+1/2}}} \,, \quad x_{j}\,\leq \,
x\,\leq\,x_{j+1}\,. $

\smallskip \noindent 
Then if $\,v\in  H^{1}(0,1)\,$ and $\,q \in H^{2}(0,1)\,$, we have the
interpolation error estimates~:

\setbox19=\hbox{$\displaystyle {{{\rm d}v}\over{{\rm d}x}}$}
\setbox20=\hbox{$\displaystyle \nor{\box19} \ib{0}$}
\smallskip \noindent   (5.24) $\qquad \displaystyle
\parallel v-\MT v \parallel \ib{0} \,\,\,\leq\,\, C\,\hT\, \box20 \, $
\setbox29=\hbox{$\displaystyle{{{\rm d}^{2}q}\over{{\rm d}x^{2}}} $}
\setbox30=\hbox{$\displaystyle \nor{\box29} \ib{0}$}

\smallskip \noindent   (5.25) $\qquad \displaystyle 
\parallel q- \PiT q \parallel \ib{1} \,\,\,\leq\,\, C\,\hT\, \box30  \,$

\smallskip \noindent 
where $\,\hT\,$,  defined in (4.13), is the maximal step size in mesh
$\,\,\cal{T}\,$ and $\,C\,$ is some constant independant of $\,\,\cal{T}\,$,
$\,v\,$ and $\,q\,$.

\bigskip
\noindent {\bf Proof of Theorem 4.}

\noindent $\bullet \quad$
First the Poisson equation (5.5) is formulated under the Petrov-Galerkin form
(2.3)-(2.4) in linear space $\,V= L^{2}(0,1)\times H^{1}(0,1)\,$. Then
Proposition~1 about continuous inf-sup condition and infinity condition and
Theorem~1 show that the first hypothesis of Theorem~2 is satisfied. 

\noindent $\bullet \quad$
Secondly let $\,{\cal{U}}_{\alpha,\beta}\,$ be a family of uniformly regular
meshes $\,\cal{T}\,$. The discrete inf-sup condition is satisfied with a
constant $\,\delta\,$ in the right hand side of (4.6) which does not depend on 
 $\,\cal{T}\,$, due to Proposition 7 and in particular inequality (5.16).
 
\noindent $\bullet \quad$
We prove now the infinity condition (4.7) between $\,V_{1}= \UT\times\PT\,$ and 
 $\,V_{2}= \UT\times\QT^{\psi}\,$. Let $\,\eta=(v,q)\,$ be a non-zero pair in 
$\,V_{2}\,$.

 $\star \quad$ If $\,{\rm div}\,q\ne 0\,$, let $\,u=\lambda\,{\rm div}\,q\,$ and 
 $\,p=0 . $ We set $\,\xi=(u,p)\in \UT\times\PT\,$ and we have
$\,\gamma(\xi,\eta)=\lambda\,\parallel {\rm div}\,q \parallel^{2}\,$ which tends to
$+\infty\,$ when $\,\lambda\,$ tends to infinity. 

 $\star \quad$ If $\,{\rm div}\,q =  0\,$, and $\,v\ne0\,$, we construct $\,p\,$ as
the linear interpolate of $\,{\rm grad}\,\varphi\,$, where $\,\varphi\in
H^{1}_{0}(0,1)\,$ is the variational solution of Poisson problem $\,\Delta
\varphi=v\,$. Then $\displaystyle \,(p,q) = \Bigl( \int_{0}^{1}p(x)\,{\rm d}x
\Bigr)\,q\,\, $  because  $\,{\rm div}\,q =  0\,$ implies that $\,q\,$ is equal
to some constant.  But
 
\setbox18=\hbox{$\displaystyle \int_{0}^{1}p(x)\,{\rm d}x  $}
\setbox20=\hbox{$\displaystyle \sum_{j=1}^{n-1} \int_{\displaystyle
x_{j}}^{\displaystyle x_{j+1}} \biggl( {{{\rm d}\varphi}\over {{\rm
d}x}}(x_{j})\, {{x_{j+1}-x}\over {h_{j+1/2}}} +  {{{\rm d}\varphi}\over {{\rm
d}x}}(x_{j+1})\, {{x-x_{j}}\over {h_{j+1/2}}}\, \biggr)\, {\rm d}x $}
\setbox21=\hbox{$\displaystyle \sum_{j=1}^{n-1} \int_{\displaystyle
x_{j}}^{\displaystyle x_{j+1}} {{{\rm d}\varphi}\over {{\rm d}x}}(x)  {\rm d}x$}
\setbox37=\vbox{\halign{#&#&#\cr \box18&$\,\,=\,\,$&\box20 \cr 
&$\,\,=\,\,$& \box21 \cr}}

\smallskip \noindent  $\box37$

\smallskip \noindent
because ${\rm grad}\,\varphi$ is affine in each element $\,]x_{j},x_{j+1}[\,$
since $\,\Delta\,\varphi=v\,$ is a constant in each such interval. 
We deduce that  $\displaystyle \,\int_{0}^{1}p(x)\,{\rm d}x=0\,$ due to the
homogeneous Dirichlet boundary conditions for function $\,\varphi\,$. We take
$\,\xi=(0 \,,\, \lambda\,p)\,$. Then 

\smallskip \noindent   $ \displaystyle 
\gamma(\xi\,,\,\eta)\,\,=\,\, (\lambda\,p\,,\,q) + (0\,,\,{\rm div}\,q) + (
\lambda\,{\rm div}\,p\,,\,v)\,\,=\,\, \lambda \parallel v \parallel \ib{0}^{2} \,$

\smallskip \noindent
and this expression tends towards $\,+\infty\,$ as  $\,\lambda\,$ tends to 
$\,+\infty\,$. 

 $\star \quad$ If $\,{\rm div}\,q\,$ and $\,v\,$ are both equal to zero, $\,q\,$ is a
constant function which is not null because $\,\eta \ne0\,$. If we take
$\,u=0\,$ and $\,p=\lambda\,q\,$ (this last choice is possible because, due to
(5.1), $\,\PT\,$ and $\,\QT^{\psi}\,$  contain the constant functions), we get 
$\,\gamma(\xi,\eta)\,\,=\,\,  \lambda \parallel q \parallel \ib{0}^{2}\,$ and
this expression tends to  $\,+\infty\,$  as  $\,\lambda\,$ tends to 
$\,+\infty\,$. Therefore the discrete infinity condition (4.7) is satisfied. 

\smallskip \noindent $\bullet \quad$
The conclusion of Theorem 2 ensures the majoration of the error in 
$\, L^{2}(0,1)\times H^{1}(0,1)\,$ norm (left hand side of relations (4.8)
and (5.7)) by the interpolation error (right hand side of relation (4.8)). From
Proposition~8, the interpolation error is of order one and we have 

\setbox19=\hbox{$\displaystyle  {{{\rm d}u}\over{{\rm d}x}}$}
\setbox20=\hbox{$\displaystyle  {{{\rm d}^{2}p}\over{{\rm d}x^{2}}}$}
\smallskip \noindent   (5.26) $\qquad \displaystyle 
\parallel u-\uT \parallel \ib{0} + \parallel p-\pT \parallel \ib{1} 
\quad\leq\quad C \,\, \hT\,\, \biggl( \,\nor{\box19}_{0} +\nor{\box20}_{0}
\, \biggr)  \,$

\smallskip \noindent
when  $\,\cal{T}\,$ belongs to family   $\,{\cal{U}}_{\alpha,\beta}\,$ of
uniformly regular meshes. The final estimale (5.7) is a consequence of regularity
of the solution $\,u\,$ of the homogeneous Dirichlet Poisson problem (5.5) when
$\,f\,$ belongs to $\,H^{1}(0,1)\,$ : 

\smallskip \noindent   (5.27) $\qquad \displaystyle 
u \in H^{3}(0,1)\quad{\rm and}\quad \parallel u \parallel \ib{3}
\,\,\leq\,\, \widetilde{C} \, \parallel f \parallel \ib{1}  \,$

\smallskip \noindent   
Joined with (5.26), this inequality ends the proof of Theorem~4.
 $\hfill  \square \kern0.1mm$

\bigskip \bigskip   \noindent {\bf \large 6) \quad  First order for least squares} 

\monitem  
We have established with Theorem 4  that convergence of the finite
volume method but the result suffers from the fact that a too important
regularity is necessary for the datum of homogneous Dirichlet problem of Poisson
equation

\smallskip \noindent   (6.1) $\qquad \displaystyle 
-\Delta \,u\,\,=\,\,f \quad {\rm in}\,\,]0,1[\,, \qquad u(0)\,=\,u(1)\,=\,0\,. $

\smallskip \noindent 
The dream would be to use the interpolation result

\smallskip \noindent   (6.2) $\qquad \displaystyle 
\setbox19=\hbox{$\displaystyle  {{{\rm d}u}\over{{\rm d}x}}$}
\parallel u - \MT u \parallel\ib{0}\quad\leq \quad C\,\,\hT\,\,\nor{\box19}\ib{0}
\,$

\smallskip \noindent 
but if $\,u\,$ belongs only in $\,H^{1}_{0}(0,1) ,\,$ its gradient
$\displaystyle \,\,p={{{\rm d}u}\over{{\rm d}x}}\,\,$ belongs only in $\,L^{2}(0,1)\,$ and
there is no hope to define the interpolate $\,\PiT p\,$ for a so poor regular
function and consequently to define fluxes at interfaces between two finite
elements or the   $\,L^{2}(0,1)\,$ scalar product $\,(f,v)\,$. 

\smallskip \noindent $\bullet \quad$ 
Secondly, the finite element method with linear finite elements show both
estimates [CR72] :

\setbox19=\hbox{$\displaystyle {{{\rm d}^{2}u}\over{{\rm d}x^{2}}}$}

\smallskip \noindent   (6.3) $\qquad \displaystyle 
\parallel u - \uT \parallel\ib{1}\quad\leq \quad C\,\,\hT\,\,\nor{\box19}\ib{0} $
\setbox29=\hbox{$\displaystyle  {{{\rm d}^{2}u}\over{{\rm d}x^{2}}}$}

\smallskip \noindent   (6.4) $\qquad \displaystyle 
\parallel u - \uT \parallel\ib{0}\quad\leq \quad C\,\,\hT^{2}
\,\,\nor{\box29}\ib{0} \,. $

\smallskip \noindent
Inequality (6.3) is not accessible for present finite volumes because the
discrete unknown field $\,\uT\,$ belongs only in $\,L^{2}(0,1)\,$ and estimate
(6.4) show second order accuracy in the $\,L^{2}\,$  norm, which is much more
precise than the interpolation estimate (6.2) can do. We will show in next
theorem that the intermediate result 

\smallskip \noindent  $ \displaystyle 
\parallel u - \uT \parallel\ib{0}\quad\leq \quad C\,\,\hT \,\,\parallel f
\parallel\ib{0} \,$

\smallskip \noindent 
holds when $\,f\,$ belongs  in  $\,L^{2}(0,1)\,$. This result is optimal in
the sense that on one hand the   $\,H^{2}\,$ semi-norm in the right hand side of
(6.3) and (6.4) demands a minimum of regularity for  datum $\,f\,$ and
condition $\,f\in  L^{2}(0,1)\,$ is a good regularity constraint for a
distribution which {\it a priori} belongs to space  $\,H^{-1}(0,1)\,$. On the
other hand, the $\,L^{2}\,$ error $\,\parallel u - \uT \parallel \ib{0}\,$ should
have the same order that the interpolation error $\,\parallel u - \MT u \parallel
\ib{0}\,$ (see left hand side of (6.2)).

\smallskip \noindent $\bullet \quad$ 
Nevertheless, note that some kind of superconvergence between the interpolated
value $\, \MT u \,$ and the discrete solution $\, \uT ,\,$ {\rm i.e.} estimation of
the type 

\smallskip \noindent  $ \displaystyle 
\parallel \MT u - \uT \parallel\ib{0}\quad\leq \quad C\,\,\hT^2 \,$

\smallskip \noindent
have been obtained by Arbogast, Wheeler and Yotov [AWY97] in the case of
quasi-uniform grids and sufficiently regular solution $\,u .\,$

\bigskip
\noindent {\bf Theorem 5.  $\quad$ A second result of convergence.}

\noindent 
We make the same hypotheses than in Theorem 4 for the interpolation function
$\,\psi\,$, for the family  $\cal{U}_{\alpha,\beta}$ $(0<\alpha<1<\beta\,)$ of
uniformly regular meshes $\cal T$ and we suppose that datum $\,f\in L^{2}(0,1) \,$
is  given. Then the solution $\,u\in H^{2}(0,1) \,$ of problem (6.1) can be
approximated by the finite volume method 

\setbox19=\hbox{$\displaystyle  \quad \xiT = (\uT,\pT) \in \UT\times \PT$}
\setbox20=\hbox{$\displaystyle  \quad \gamma(\xiT,\eta) \,\,=\,\, (f,v)\,, \quad
\forall \, \eta \,\in \UT\times\QT^{\psi}$}

\smallskip \noindent   (6.5) $\qquad \displaystyle 
\left\{ \begin{array} {l}  \displaystyle  
 \box19 \\ \box20 
 \end{array} \right.  $ 

\smallskip \noindent 
with $\,\UT\,,\PT\,,\QT^{\psi}\,$ and
$\,\gamma(\smb,\smb)\,$  defined  in (1.11), (1.12), (3.6) and (2.5)
respectively.  Moreover there exists some constant
$\,C\,$ depending only on $\,\alpha\,$ and $\,\beta\,$ such that 

\smallskip \noindent   (6.6) $\qquad \displaystyle 
\parallel u - \uT \parallel\ib{0}+ \parallel p - \pT \parallel\ib{0}\quad\leq
\quad C\,\,\hT \,\,\parallel f \parallel\ib{0} \,, $

\smallskip \noindent
 with $\,\hT\,$ equal to the maximal  size of mesh $\,\cal{T}\,$.

\bigskip
\noindent {\bf Proposition  9.  $\quad$ Complementary interpolation estimate.}

\noindent 
Let $\,q\,$ be a given function in $\,H^{1}(0,1)\,$ and $\,\PiT q\,$ be its
linear interpolate in space $\,\PT\,$ associated with the mesh  $\,\cal{T}\,$
and defined in (5.23). Then we have

\smallskip \noindent   (6.7) $\qquad \displaystyle 
\setbox20=\hbox{$\displaystyle {{{\rm d}q}\over{{\rm d}x}}$}
\parallel q - \PiT q \parallel \ib{0} \quad \leq \quad C \,\, \hT\,\,
\nor{\box20} \ib{0} \,$

\smallskip \noindent
where $\,\hT\,$ is  the maximal  step size of mesh $\,\cal{T}\,$ and $\,C\,$ some
constant independent of  $\,\cal{T}\,$ and $\,q\,$. 

\bigskip \noindent  {\bf Proof of proposition 9.}

\monitem 
The proof of this proposition is conducted as in Proposition 8. We first
establish inequality (6.7) when $\,{\cal{T}}=\{0=x_{0}<x_{1}=1 \}\,$ is the
trivial mesh of interval $\,]0,1[\,$. In this particular case, function $\,\,q  -
\PiT q\,\,$ belongs to $\,H^{1}_{0}(0,1)\,$ and the Poincar\'e estimate show that
we have

\setbox20=\hbox{$\displaystyle {{{\rm d}}\over{{\rm d}x}} \bigl(  q - \PiT q
\bigr) $}

\smallskip \noindent   (6.8) $\qquad \displaystyle 
\parallel q - \PiT q \parallel \ib{0} \quad \leq \quad C_{1} \,\,
\nor{\box20}\ib{0}  \,\,. \,$

\smallskip \noindent
Then we can establish the simple estimation

\setbox19=\hbox{$\displaystyle {{{\rm d}}\over{{\rm d}x}} \bigl(  \PiT q \bigr)$}
\setbox20=\hbox{$\displaystyle {{{\rm d}q}\over{{\rm d}x}} $}
\smallskip \noindent   (6.9) $\qquad \displaystyle 
\nor{\box19}\ib{0}   \quad \leq \quad  \,\, \nor{\box20}\ib{0} \,$

\smallskip \noindent
because

\setbox19=\hbox{$\displaystyle {{{\rm d}}\over{{\rm d}x}} \bigl( \PiT q \bigr) $}
\setbox20=\hbox{$\displaystyle \nor{\box19}_{0}^{2}$}
\setbox21=\hbox{$\displaystyle \int_{0}^{1} \bigl( q(1)-q(0) \bigr)^{2}\,{\rm d}x$}
\setbox22=\hbox{$\displaystyle \biggl( \int_{0}^{1} {{{\rm d}q}\over{{\rm d}y}}\,
{\rm d} y \, \biggr)^{2} $}
\setbox23=\hbox{$\displaystyle \int_{0}^{1} \Bigl( {{{\rm d}q}\over{{\rm d}y}}
\Bigr)^{2}\,{\rm d} y$}
\setbox25=\hbox{$\displaystyle {{{\rm d}q}\over{{\rm d}x}} $}
\setbox24=\hbox{$\displaystyle \nor{\box25} \ib{0}^{2}\,.$}
\setbox37=\vbox{\halign{#&#&#&#&#\cr 
\box20 &$\quad=\quad$& \box21 &$\quad=\quad$ &\box22 \cr 
  &$\quad\leq\quad$& \box23 &$\quad=\quad$ &\box24 \cr }}

\smallskip \noindent $\box37$

\smallskip \noindent
The proof of estimate (6.7) in this particular case follows from triangular
inequality based on (6.8) and (6.9) with $\,C=2\,C_{1}\,$.

\bigskip \noindent $\bullet \quad$
A general mesh $\,{\cal{T}} = \{0=x_{0}<x_{1}<\cdots<x_{n}=1 \}\,$ is composed
with $\,n\,$ trivial meshes $\,{\cal{T}}_{j+1/2} = \{ x_{j}<x_{j+1}\}\,$
of the interval $\,]x_{j},x_{j+1}[\,$. We adopt the notation (5.34)
introduced inside the proof of Proposition 8 and we have :

\setbox19=\hbox{$\displaystyle \parallel q - \PiT q \parallel \ib{0,]0,1[}^{2}$}
\setbox20=\hbox{$\displaystyle q - \PiT q $}
\setbox21=\hbox{$\displaystyle \sum _{j=0}^{n-1} \, \nor{\box20}
\ib{0,]x_{j},x_{j+1}[}^{2}$}
\setbox22=\hbox{$\displaystyle \widehat{q}_{j+1/2} -\widehat{\Pi}
\widehat{q}_{j+1/2} $}
\setbox23=\hbox{$\displaystyle   \sum _{j=0}^{n-1}  h_{j+1/2} \, \nor{\box22} 
\ib{0,]0,1[}^{2} $}
\setbox24=\hbox{$\,\,$ from (5.36)}
\setbox25=\hbox{$\displaystyle {{{\rm d}}\over{{\rm d}\theta}}\, \Bigl( \,
\widehat{q}_{j+1/2} -\widehat{\Pi} \widehat{q}_{j+1/2}\, \Bigr) $}
\setbox26=\hbox{$\displaystyle   (C_{1})^{2}\, \sum _{j=0}^{n-1}  h_{j+1/2} \, 
\nor{\box25}  \ib{0,]0,1[}^{2} $}
\setbox27=\hbox{$\,\,$ from (6.8)}
\setbox28=\hbox{$\displaystyle {{{\rm d}}\over{{\rm d}x}}\, \bigl( \, q - \PiT q \bigr) $}
\setbox29=\hbox{$\displaystyle   (C_{1})^{2}\, \sum_{j=0}^{n-1}  h_{j+1/2}^{2}
\,  \nor{\box28}  \ib{0,]x_{j},x_{j+1}[}^{2} $}
\setbox30=\hbox{$\,\,$ from (5.37).}
\setbox37=\vbox{\halign{#&#&#&#\cr  \box19 &$\,\,=\,\,$& \box21& \cr 
 &$\,\,=\,\,$& \box23  & \box24 \cr  
 &$\,\,\leq\,\,$& \box26  & \box27 \cr 
  }}

\smallskip \noindent $\box37$

\setbox19=\hbox{$\displaystyle \parallel q - \PiT q \parallel \ib{0,]0,1[}^{2}$}
\setbox20=\hbox{$\displaystyle q - \PiT q $}
\setbox21=\hbox{$\displaystyle \sum _{j=0}^{n-1} \, \nor{\box20}
\ib{0,]x_{j},x_{j+1}[}^{2}$}
\setbox22=\hbox{$\displaystyle \widehat{q}_{j+1/2} -\widehat{\Pi}
\widehat{q}_{j+1/2} $}
\setbox23=\hbox{$\displaystyle   \sum _{j=0}^{n-1}  h_{j+1/2} \, \nor{\box22} 
\ib{0,]0,1[}^{2} $}
\setbox24=\hbox{$\,\,$ from (5.36)}
\setbox25=\hbox{$\displaystyle {{{\rm d}}\over{{\rm d}\theta}}\, \Bigl( \,
\widehat{q}_{j+1/2} -\widehat{\Pi} \widehat{q}_{j+1/2}\, \Bigr) $}
\setbox26=\hbox{$\displaystyle   (C_{1})^{2}\, \sum _{j=0}^{n-1}  h_{j+1/2} \, 
\nor{\box25}  \ib{0,]0,1[}^{2} $}
\setbox27=\hbox{$\,\,$ from (6.8)}
\setbox28=\hbox{$\displaystyle {{{\rm d}}\over{{\rm d}x}}\, \bigl( \, q - \PiT q \bigr) $}
\setbox29=\hbox{$\displaystyle   (C_{1})^{2}\, \sum_{j=0}^{n-1}  h_{j+1/2}^{2}
\,  \nor{\box28}  \ib{0,]x_{j},x_{j+1}[}^{2} $}
\setbox30=\hbox{$\,\,$ from (5.37).}
\setbox37=\vbox{\halign{#&#&#&#\cr  \box19 &$\,\,\leq\,\,$& \box29  & \box30 \cr }}

\smallskip \noindent $\box37$

\smallskip \noindent 
Then

\setbox28=\hbox{$\displaystyle {{{\rm d}}\over{{\rm d}x}}\, \bigl( \, q - \PiT q
\bigr) $}
\smallskip \noindent   (6.10) $\qquad \displaystyle 
\parallel q - \PiT q \parallel \ib{0} \quad\leq\quad C_{1}\,\,\hT\,\,\nor{\box28}
\ib{0}  \, . $

\smallskip \noindent  
In an analogous way than the one that conducted to estimation (6.9), we have~:
\setbox19=\hbox{$\displaystyle {{{\rm d}}\over{{\rm d}x}}  \bigl( \PiT q \bigr) $}
\setbox20=\hbox{$\displaystyle \nor{\box19}\ib{0}^{2}$}
\setbox21=\hbox{$\displaystyle \sum_{j=0}^{n-1} \int_{\displaystyle
x_{j}}^{\displaystyle x_{j+1}} \biggl( {{q(x_{j+1})-q(x_{j})}\over{h_{j+1/2}}}
\biggr) ^{2}\,{\rm d}x$}
\setbox22=\hbox{$\displaystyle \sum_{j=0}^{n-1} {{1}\over{h_{j+1/2}}} \,\biggl( 
\int_{\displaystyle x_{j}}^{\displaystyle x_{j+1}}  \biggl( {{{\rm d}q}\over{{\rm
d}x}}\biggr) \,{\rm d}x  \biggr) ^{2} $}
\setbox23=\hbox{$\displaystyle \sum_{j=0}^{n-1} {{1}\over{h_{j+1/2}}} 
\biggl( \int_{\displaystyle x_{j}}^{\displaystyle x_{j+1}}  \biggl( {{{\rm
d}q}\over{{\rm d}x}}\biggr) ^{2} \,{\rm d}x  \biggr) 
\, \Bigl( \int_{\displaystyle x_{j}}^{\displaystyle x_{j+1}}  \,{\rm d}x \Bigr)$}
\setbox24=\hbox{$\quad$ }
\setbox25=\hbox{$\displaystyle  {{{\rm d}q}\over{{\rm d}x}}$}
\setbox26=\hbox{$\displaystyle \nor{\box25}\ib{0}^{2}$}
\setbox37=\vbox{\halign{#&#&#&#\cr 
\box20 &$\,\,=\,\,$& \box21 &  \cr       &$\,\,=\,\,$ &\box22 & \cr 
  &$\,\,\leq\,\,$& \box23  &\box24 \cr   & & & by Cauchy-Schwarz  \cr  
   &$\,\,=\,\,$& \box26  & \cr  }}

\smallskip \noindent $\box37$ 

\setbox19=\hbox{$\displaystyle {{{\rm d}}\over{{\rm d}x}}  \PiT q $}
\setbox25=\hbox{$\displaystyle  {{{\rm d}q}\over{{\rm d}x}}$}
\smallskip \noindent   (6.11) $\qquad \displaystyle 
\nor{\box19}\ib{0} \quad\leq\quad \nor{\box25}\ib{0}\,.\,$

\smallskip \noindent  
Then inequality (6.10) joined with (6.11) and the
triangular inequality show (6.7) with $\,C=2\,C_{1}\,$.
$\hfill  \square \kern0.1mm$

\bigskip
\noindent {\bf Proof of Theorem 5.}

\noindent $\bullet \quad$
We divide it into three steps. First we establish that if a pair $\,(\sT,\mT)
\in \UT\times\PT\,$ is solution of the discrete finite volume problem in
Petrov-Galerkin formulation, with data $\,\delta\,$ and $\,\varphi\,$ in
$\,L^{2}(0,1)\,$

\smallskip \noindent   (6.12) $\qquad \displaystyle 
(\mT,q) + (\sT,{\rm div}\,q)\,\,=\,\, (\delta,q) + (\varphi,{\rm div}\,q)\,,\quad
\forall  \,q\in \QT^{\psi} \,$

\smallskip \noindent   (6.13) $\qquad \displaystyle 
({\rm div}\,\mT,v)\,\, \qquad \qquad =\,\,0\,,\qquad \qquad \qquad  \quad \,\,\,\,
\forall  \,v\in \UT \,$

\smallskip \noindent 
then we have a stability estimate

\smallskip \noindent   (6.14) $\qquad \displaystyle 
\parallel \sT \parallel \ib{0} + \parallel \mT \parallel \ib{1} \quad
\leq \quad C \,\Bigl( \parallel \delta  \parallel \ib{0} + \parallel \varphi 
\parallel \ib{0} \Bigr) \,$

\smallskip \noindent  
where $\,C\,$ is a constant dependent only on parameters $\,\alpha\,,\beta\,$ of
the class  $\,\cal{U}_{\alpha,\beta}\,$ of uniform meshes. Since $\,\psi\,$
interpolant function satisfies the interpolation compatibiliy condition,
Proposition 7 establishes that the discrete inf-sup condition is uniformly
satisfied : 
 
\setbox21=\hbox{$\quad \exists \,\rho\,>0\,,\quad \forall\,{\cal{T}}\!\!\in
{\cal{U}}_{\alpha,\beta}\,,\quad\forall\, \xi \, =\, (u,p) \in
\UT\times\PT\,,   \xi \ne 0\,,  $} 
\setbox22=\hbox{$\quad \exists\,\eta=(v,q)\in\UT\times\QT^{\psi}\,,\quad 
\parallel  \eta \parallel\,\leq\,1 \quad {\rm and}\quad
\gamma(\xi,\eta)\,\,\geq\,\,\rho\,\parallel  \xi \parallel\, .$}

\smallskip \noindent   (6.15) $\qquad \displaystyle 
\left\{ \begin{array} {l}  \displaystyle  
 \box21 \\ \box22 
 \end{array} \right.  $ 

\smallskip \noindent  
We use this stability inequality with $\,\xi=(\sT,\mT)\,$ solution of problem
(6.12)-(6.13). Then there exists $\,\eta=(v,q)\in \UT\times\QT^{\psi}\,$ such
that $\,\parallel \eta \parallel \,\leq\,1\,$ and

\setbox19=\hbox{$\displaystyle  {{1}\over{\sqrt{2}}} \, \bigl( \parallel \sT
\parallel \ib{0} + 
\parallel \mT \parallel \ib{1} \bigr) \,\, \leq \,\,\,\,  \parallel \xi
\parallel$}
\setbox21=\hbox{$\displaystyle {{1}\over{\rho}}\,\gamma(\xi,\eta)$}
\setbox22=\hbox{$\displaystyle {{1}\over{\rho}}\,\Bigl( 
(\delta \,,\, q) + (\varphi \,,\, {\rm div}\,q) \Bigr) $}
\setbox23=\hbox{$\displaystyle {{1}\over{\rho}}\,\bigl(
\parallel \delta \parallel \ib{0} \, \parallel q \parallel \ib{0} + \parallel
\varphi \parallel\ib{0}  \parallel {\rm div}\,q \parallel\ib{0} \bigr) $}
\setbox24=\hbox{$\displaystyle {{1}\over{\rho}}\,\bigl(\parallel \delta \parallel
\ib{0}  + \parallel \varphi \parallel \ib{0} \bigr)  \parallel \eta \parallel $}
\setbox37=\vbox{\halign{#&#&#\cr 
\box19 &$ \quad\leq\quad $& \box21   \cr       &$ \quad=\quad$& \box22   \cr 
  &$ \quad\leq\quad$& \box23 \cr    }}
\smallskip \noindent  $\box37$

\setbox19=\hbox{$\displaystyle  {{1}\over{\sqrt{2}}} \, \bigl( \parallel \sT
\parallel \ib{0} +  \parallel \mT \parallel \ib{1} \bigr) \,\, $}
\setbox21=\hbox{$\displaystyle {{1}\over{\rho}}\,\gamma(\xi,\eta)$}
\setbox22=\hbox{$\displaystyle {{1}\over{\rho}}\,\Bigl( 
(\delta \,,\, q) + (\varphi \,,\, {\rm div}\,q) \Bigr) $}
\setbox23=\hbox{$\displaystyle {{1}\over{\rho}}\,\bigl(
\parallel \delta \parallel \ib{0} \, \parallel q \parallel \ib{0} + \parallel
\varphi \parallel\ib{0}  \parallel {\rm div}\,q \parallel\ib{0} \bigr) $}
\setbox24=\hbox{$\displaystyle {{1}\over{\rho}}\,\bigl(\parallel \delta \parallel
\ib{0}  + \parallel \varphi \parallel \ib{0} \bigr)  \parallel \eta \parallel $}
\setbox37=\vbox{\halign{#&#&#\cr 
\box19 &$ \quad\leq\quad $  & \box24   \cr  }}
\smallskip \noindent  $\box37$

\smallskip \noindent  
and inequality (6.14) is a direct consequence of the fact that $\,\parallel \eta
\parallel \,\leq\,1\,.$

\bigskip \noindent $\bullet \quad$
Secondly let $\,\wT\,$ and $\,\muT\,$ be two arbitrary functions  in spaces
$\,\UT\,$ and $\,\PT\,$ respectively. From the continuous mixed formulation 

\smallskip \noindent   (6.16) $\qquad \displaystyle 
(p\,,\,q) + (u\,,\,{\rm div}\,q)\,\,=\,\,0\,\quad \forall\, q \in H^{1}(0,1) \,$

\smallskip \noindent   (6.17) $\qquad \displaystyle 
({\rm div}\,p \,,\, v) + (f\,,\,v)\,\,=\,\,0\,\quad \forall\, v \in L^{2}(0,1) \,$

\smallskip \noindent 
and the discrete Petrov-Galerkin approximation

\smallskip \noindent   (6.18) $\qquad \displaystyle 
(\pT \,,\,q) + (\uT \,,\,{\rm div}\,q)\,\,=\,\,0\,\quad \forall\, q \in \QT^{\psi} \,$ 

\smallskip \noindent   (6.19) $\qquad \displaystyle 
({\rm div} \, \pT \,,\,v) + (f \,,\,v)\,\,=\,\,0\,\quad \forall\, v \in \UT \,. $

\smallskip \noindent 
We deduce by difference

\smallskip \noindent   (6.20) $\quad \displaystyle 
(\pT - \muT \,,\, q) + (\uT - \wT  \,,\,{\rm div}\,q)\,=\,(p - \muT \,,\,q) +
(u -\wT  \,,\,{\rm div}\,q) \,\,\, \forall\, q \in \QT^{\psi} \,$

\smallskip \noindent   (6.21) $\qquad \displaystyle 
({\rm div} \, (\pT - \muT) \,,\,v) \,\,=\,\,({\rm div} \, (p - \muT)
\,,\,v)\,\quad \forall\, v \in \UT \,. $

\smallskip \noindent 
If we select for $\,\muT\,$ the $P_{1}$ interpolate of $\,p\,$ in space
$\,\PT\,$, {\it i.e.}  $\,\muT = \PiT p\,$, we have $\,p(x_{j})=\PiT p(x_{j})\,$
for each vertex $\,x_{j}\,$ of mesh $\,\cal{T}\,$, then 

\smallskip \noindent   $ \displaystyle 
\int_{\displaystyle x_{j}}^{\displaystyle x_{j+1}} {\rm div} \bigl( p - \PiT p
\bigr)\, {\rm d}x \quad=\quad 0 \,$

\smallskip \noindent 
and the same property is true for the right hand side of (6.21). Considering now
the particular case of $\,\wT=\MT u\,$, we deduce from (6.20)(6.21) and previous
estimate (6.14) the inequality

\smallskip \noindent   (6.22) $\quad \displaystyle 
\parallel \pT - \PiT p \parallel \ib{1} + \parallel \uT - \MT u \parallel
\ib{0} \quad \leq \quad C \, \Bigl(   \parallel p - \PiT p \parallel \ib{0} +
\parallel u - \MT u \parallel \ib{0}  \bigr) \,. $

\smallskip \noindent
Joined with the triangular inequality and majoration of $\,L^{2}\,$ norm by the
$\,H^{1}\,$ norm, we obtain

\smallskip \noindent   (6.23) $\qquad \displaystyle 
\parallel p - \pT  \parallel \ib{0} + \parallel u - \uT \parallel
\ib{0} \quad \leq \quad (1+C) \,  \Bigl(   \parallel p - \PiT p \parallel \ib{0} +
\parallel u - \MT u \parallel \ib{0}  \bigr) \,. $

\bigskip \noindent $\bullet \quad$
The end of the proof is a direct consequence of Propositions 8 and 9 and in
particular estimations (5.24) and (6.7) :

\setbox19=\hbox{$\displaystyle {{{\rm d}p}\over{{\rm d}x}} $}
\smallskip \noindent   (6.24) $\qquad \displaystyle 
\parallel p - \pT  \parallel \ib{0} + \parallel u - \uT \parallel
\ib{0} \quad \leq \quad C\,\,\hT\,\,   \Bigl( \,  \nor{\box19}
 \ib{0} + \parallel p \parallel \ib{0} \, \bigr) \,$
 
\smallskip \noindent
joined with the classical estimate that comes from the variational formulation
of  problem (6.1) :

\smallskip \noindent   (6.25) $\qquad \displaystyle 
\parallel p \parallel \ib{1} \quad \leq \quad C \,\, \parallel f \parallel
\ib{0} \,. $

\smallskip \noindent 
The sequence of inequalities  (6.24) and (6.25) establishes completely the
inequality (6.6) modulo classical conventions in numerical analysis concerning
the so-called {\it constant} C.    $\hfill  \square \kern0.1mm$

\bigskip \bigskip  
\noindent {\bf \large 7) \quad Conclusion and aknowledgments }  

\monitem  
In two space dimensions, mass lumping of mass matrix of mixed finite elements
has defined a particular finite volume method analysed by Baranger et al
[BMO96]. Note also that  first results for the Laplace equation approached by
a simple finite volume method on  Delaunay-Vorono\"{\i}  triangular meshes have
been obtained  by   Herbin [He95].  We think also that our one-dimensional result
for finite volumes via mixed Petrov-Galerkin finite elements  can be generalized
in dimension 2 and 3 for regular triangular or tetrahedral meshes with a numerical
scheme like the diamond scheme suggested by Noh many years ago [No64], first
analysed by Coudi\`ere, Vila and Villedieu [CVV99],  or our ``wedding scheme''
proposed in an other context~[Du92].

\bigskip \noindent $\bullet \quad$ 
We have been introduced to techniques of Petrov-Galerkin formulations thanks to
a pedagogical initiative of Bernard Larrouturou at Ecole Polytechnique. The
breakthrough of this research was done during a spring school at Les Houches in
may 1996 ; we thank all the participants and in particular Olga Cueto for good
working sollicitation. This article has been  typed with $\,\TeX\,$ by the author
and we are  redevable to the competences of Jean Louis Loday. Second version of
this report is due to  particular encouragements of  Jean-Pierre Croisille. The
author thanks also the referee for helpfull suggestions.

\bigskip \bigskip   \noindent {\bf \large 8) \quad  	References }  

\smallskip

\smallskip \hangindent=9mm \hangafter=1 \noindent [Ad75] $\,\,$
R.A. Adams. {\it Sobolev spaces.} Academic Press, New York, 1975.

\smallskip \hangindent=9mm \hangafter=1 \noindent [AWY97] $\,\,$
 T. Arbogast, M.F. Wheeler, I. Yotov. Mixed finite elements for elliptic 
problems with tensor coefficients as
cell-centered finite differences, {\it SIAM J. Numer. Anal.},  vol~34, 
p.~828-852, 1997. 

\smallskip \hangindent=9mm \hangafter=1 \noindent [Ba71] $\,\,$
 I. Babu\u ska. Error-Bounds for Finite
Element Method,   {\it Numer. Math.}, vol~16, p.~322-333, 1971.

\smallskip \hangindent=9mm \hangafter=1 \noindent [BMO96] $\,\,$
 J. Baranger, J. F. Ma\^{\i}tre, F. Oudin. Connection between finite 
volumes and mixed finite element
methods, {\it Math. Mod. and Numer. Anal.}, vol~30,  p.~445-465, 1996. 

\smallskip \hangindent=9mm \hangafter=1 \noindent [Ce64] $\,\,$
  J. Cea. Approximation variationnelle
des probl\`emes aux limites,   {\it Ann. Inst. Fourier} (Grenoble), vol 14, p.
345-444, 1964. 

\smallskip \hangindent=9mm \hangafter=1 \noindent [CR72] $\,\,$
  P.G. Ciarlet, P.A. Raviart. General
Lagrange and Hermite Interpolation in~$\R^{n}$ with Applications to Finite
Element Methods,  {\it Arch.  Rational Mech. Anal.}, vol 46, p.~177-199, 1972.

\smallskip \hangindent=9mm \hangafter=1 \noindent [CVV99] $\,\,$
Y. Coudi\`ere, J.P. Vila, P. Villedieu. Convergence rate of a finite  volume scheme for a
two-dimensional convection-diffusion  problem, {\it Mod\'elisation
Math\'ematique et Analyse Num\'erique}, vol 33, p.~494-516, 1999.

\smallskip \hangindent=9mm \hangafter=1 \noindent [Du92] $\,\,$
  F. Dubois. Interpolation de Lagrange et
volumes finis. Une technique nouvelle pour calculer le gradient d'une fonction
sur les faces d'un maillage non structur\'e, {\it  Aerospatiale Espace $\&$
Defense, Internal report ST/S 104109}, february 1992.
See also {\it Lemmes finis pour la dynamique des gaz}, chapter 8, hal-00733937. 

\smallskip \hangindent=9mm \hangafter=1 \noindent [Du97] $\,\,$
  F. Dubois. Finite volumes and
Petrov-Galerkin finite elements. The unidimensional problem, {\it Institut
A\'erotechnique de Saint Cyr}, Report 295, Conservatoire National des Arts
et M\'etiers, october 1997.

\smallskip \hangindent=9mm \hangafter=1 \noindent [EGH2k] $\,\,$
 R. Eymard, T.  Gallou\"et, R.
Herbin. Finite Volume Methods,  {\it Handbook of Numerical Analysis} 
(Ciarlet-Lions Eds),  North Holland, Amsterdam, vol~7, p.~715-1022, 2000.

\smallskip \hangindent=9mm \hangafter=1 \noindent [Ga92] $\,\,$
 T.  Gallou\"et. An introduction to 
Finite Volume Methods, in {\it M\'ethodes de Volumes Finis}, cours CEA-EDF-INRIA,
Clamart, p.~1-85, 1992.

\smallskip \hangindent=9mm \hangafter=1 \noindent [He95] $\,\,$
 R. Herbin. An error estimate for a
finite volume scheme for a diffusion-convection problem in a triangular mesh,   
{\it Numer. Meth. for Part. Diff. Equations}, vol 11, p.~165-173, 1995.

\smallskip \hangindent=9mm \hangafter=1 \noindent [Hu78] $\,\,$
 T.J.R. Hughes. A Simple Scheme
for Developping ``Upwind'' Finite Elements,   {\it Int. J. of Numer.
Meth. in Eng.}, vol 12, p.~1359-1365, 1978.

\smallskip \hangindent=9mm \hangafter=1 \noindent [JN81] $\,\,$
 C. Johnson, U. N\"avert. An
analysis of Some Finite Element Methods for Advection-Diffusion
Problems, in {\it Analytical and Numerical Approaches to Asymptotic Problems in
Analysis} (Axelsson, Frank, Van des Sluis Eds), North Holland, Amsterdam,
p.~99-116, 1981. 

\smallskip \hangindent=9mm \hangafter=1 \noindent [No64] $\,\,$
W.F. Noh. CEL : A Time Dependent Two
Space Dimensional, Coupled Euler Lagrange Code, in {\it Methods in
Computational Physics}, vol 3, Academic Press, p.~117-179, 1964.

\smallskip \hangindent=9mm \hangafter=1 \noindent [Pa80] $\,\,$
 S.V. Patankar. {\it Numerical Heat Transfer
and Fluid Flow}, Hemisphere publishing, 1980.

\smallskip \hangindent=9mm \hangafter=1 \noindent [RT77] $\,\,$ 
P.A. Raviart, J.M. Thomas. 
A Mixed Finite Element Method for 2nd
Order Elliptic Problems, in {\it Lectures in Mathematics}, vol 606 (Dold-Eckmann
Eds), Springer-Verlag, Berlin, p.~292-315, 1977.

\smallskip \hangindent=9mm \hangafter=1 \noindent [TT99] $\,\,$ 
J.M. Thomas, D. Trujillo.  Analysis
of Finite Volumes Methods, {\it Universit\'e de Pau et des Pays de l'Adour,
Applied Mathematics Laboratory},  Internal Report  95-19, 1995, Mixed Finite 
Volume Methods, {\it  International J. for Numerical Methods in
Engineering}, vol 46, p.~1351-1366, 1999.

\end{document}